\documentclass{article}

\usepackage{arxiv}
\usepackage[utf8]{inputenc} % allow utf-8 input
\usepackage[T1]{fontenc}    % use 8-bit T1 fonts
\usepackage{hyperref}       % hyperlinks
\usepackage{url}            % simple URL typesetting
\usepackage{booktabs}       % professional-quality tables
\usepackage{amsfonts}       % blackboard math symbols
\usepackage{nicefrac}       % compact symbols for 1/2, etc.
\usepackage{microtype}      % microtypography
\usepackage{lipsum}		% Can be removed after putting your text content
\usepackage{graphicx}
\usepackage{natbib}
\usepackage{doi}
\newcommand*\samethanks[1][\value{footnote}]{\footnotemark[#1]}

\usepackage{algorithm}  
\usepackage{enumerate}
\usepackage{algorithmicx}
\usepackage{algpseudocode} 
\usepackage{tabularx}
\newtheorem{thm}{Theorem}[section]
\newtheorem{defn}[thm]{Definition}
\algrenewcommand{\algorithmiccomment}[1]{\hfill(#1)}
\usepackage{mathtools}
\usepackage{subfigure}
\usepackage{multirow}
\usepackage{amssymb}

\title{Simplex Initialization: A Survey of Techniques and Trends}

\author{ Mengyu Huang\thanks{These authors contributed equally to this work.} \\
	Department of Electronic and Computer Engineering\\
	The Hong Kong University of Science and Technology\\
	\texttt{mhuangak@ust.hk} \\
	\And
	Yuxing Zhong\samethanks\\
	Department of Electronic and Computer Engineering\\
	The Hong Kong University of Science and Technology\\
	\texttt{yzhongbc@ust.hk} \\
	\And
	Huiwen Yang\samethanks\\
	Department of Electronic and Computer Engineering\\
	The Hong Kong University of Science and Technology\\
	\texttt{hyangbr@ust.hk} \\
	\And
	Jiazheng Wang\\
	Theory Lab \\
	Huawei Hong Kong Research Centre\\
	\texttt{wang.jiazheng@huawei.com} \\
	\And
	Fan Zhang\\
	Theory Lab \\
	Huawei Hong Kong Research Centre\\
	\texttt{zhang.fan2@huawei.com} \\
	\And
	Bo Bai\\
	Theory Lab \\
	Huawei Hong Kong Research Centre\\
	\texttt{baibo8@huawei.com} \\
	\And
	Ling Shi\\
	Department of Electronic and Computer Engineering\\
	The Hong Kong University of Science and Technology\\
	\texttt{eesling@ust.hk} \\
}

\date{}

\renewcommand{\shorttitle}{\textit{arXiv} Template}

\hypersetup{
pdftitle={Simplex Initialization: A Survey of Techniques and Trends},
pdfsubject={Simplex Initialization},
pdfauthor={Mengyu Huang, Yuxing Zhong, Huiwen Yang, Jiazheng Wang, Fan Zhang, Bo Bai, Ling Shi},
pdfkeywords={Simplex Initialization,Linear Programming,Continuous Optimization, Learning, Algorithm Acceleration},
}

\begin{document}
\maketitle

\begin{abstract}
	The simplex method is one of the most fundamental technologies for solving linear programming (LP) problems and has been widely applied to different practical applications. In the past literature, how to improve and accelerate the simplex method has attracted plenty of research. One important way to achieve this goal is to find a better initialization method for the simplex. In this survey, we aim to provide an overview about the initialization methods in the primal and dual simplex, respectively. We also propose several potential future directions about how to improve the existing initialization methods with the help of advanced learning technologies.
\end{abstract}

% keywords can be removed

\section{Introduction}
Linear Programming (LP) is a minimization or maximization optimization problem with linear objective, linear equality and inequality constraints. The applications range from engineering, agriculture, transportation to food industry, manufacturing, etc. 

The foundation of LP dates back to the work proposed by~\cite{kantorovich1939mathematical} where an optimization problem of production planning and organization was studied. Later, Dantzig, who is known as the father of LP, introduced the first general framework of LP and provided the basic \textbf{primal simplex method} for solving LPs~\citep{dantzig1951maximization}
%~\cite{dantzig1951maximization}. 
Since then, LPs and the simplex method have been greatly explored and have led to a large number of extensions. Specifically, the \textbf{dual simplex method} \citep{lemke1954dual} was proposed in 1954. The primal and dual simplex are the two main streams in the simplex methods.

In 1972,
%Klee and Minty
~\cite{klee1972good} constructed a special LP example and showed that the basic simplex method has exponential time complexity in the worst case. From that moment on, more research efforts turned to finding more efficient algorithms, either by improving the basic simplex method or by proposing new algorithms from other perspectives. In 1984, the \textbf{interior point method (IPM)} for solving LPs, which is also called as Karmarkar’s algorithm, was proposed in~\cite{karmarkar1984new}. This algorithm is the first practically feasible method that can solve LPs in polynomial time and has prompted many studies on the variants of IPM for solving LPs. 

The variants of the simplex method (primal and dual) and the variants of the IPM are two main methods for solving LPs. In general, the simplex method is more competitive in practice, while the IPM has a better time complexity in theory. In this survey, we will first introduce and compare the main ideas of the two methods and then focus on the extension works for improving the efficiency of the simplex method. 

The reason for investigating the simplex method are threefold. First, the simplex method plays an important role in solving some very common and important optimization problems, e.g., (mixed) integer LP problems, etc.  Second, unlike IPMs, which are thought to have been theoretically and computationally matured, there still exist some research gaps for the simplex method. Last, most of the extension works of the simplex method include heuristics which are ad hoc in nature and could be further improved. In recent years, with the rapid development of machine learning and deep learning, it is interesting and meaningful to explore whether the existing heuristics can be improved by combining some learning methods, or to propose some learning-based simplex methods.

Basically, the computational time of the simplex method is mainly spent in three stages, namely, \textbf{the initialization stage, the iteration stage, and the termination stage}. In the initialization stage, an initial point or an initial basis is provided for the simplex algorithm. Then in the iteration stage, given a pivot rule, the algorithm performs the selection of the entering variable and the leaving variable to obtain the improved point or basis iteratively. Finally, in the termination stage, the algorithm is ended when some designed termination conditions are satisfied. Accordingly, the extension works of the simplex method, which to improve the algorithm efficiency and to reduce the computational cost, can be divided into three types based on which stage (initialization, iteration, or termination) is modified and improved. 

\begin{figure}
	[htbp]
	\centering
	\includegraphics[width=\linewidth]{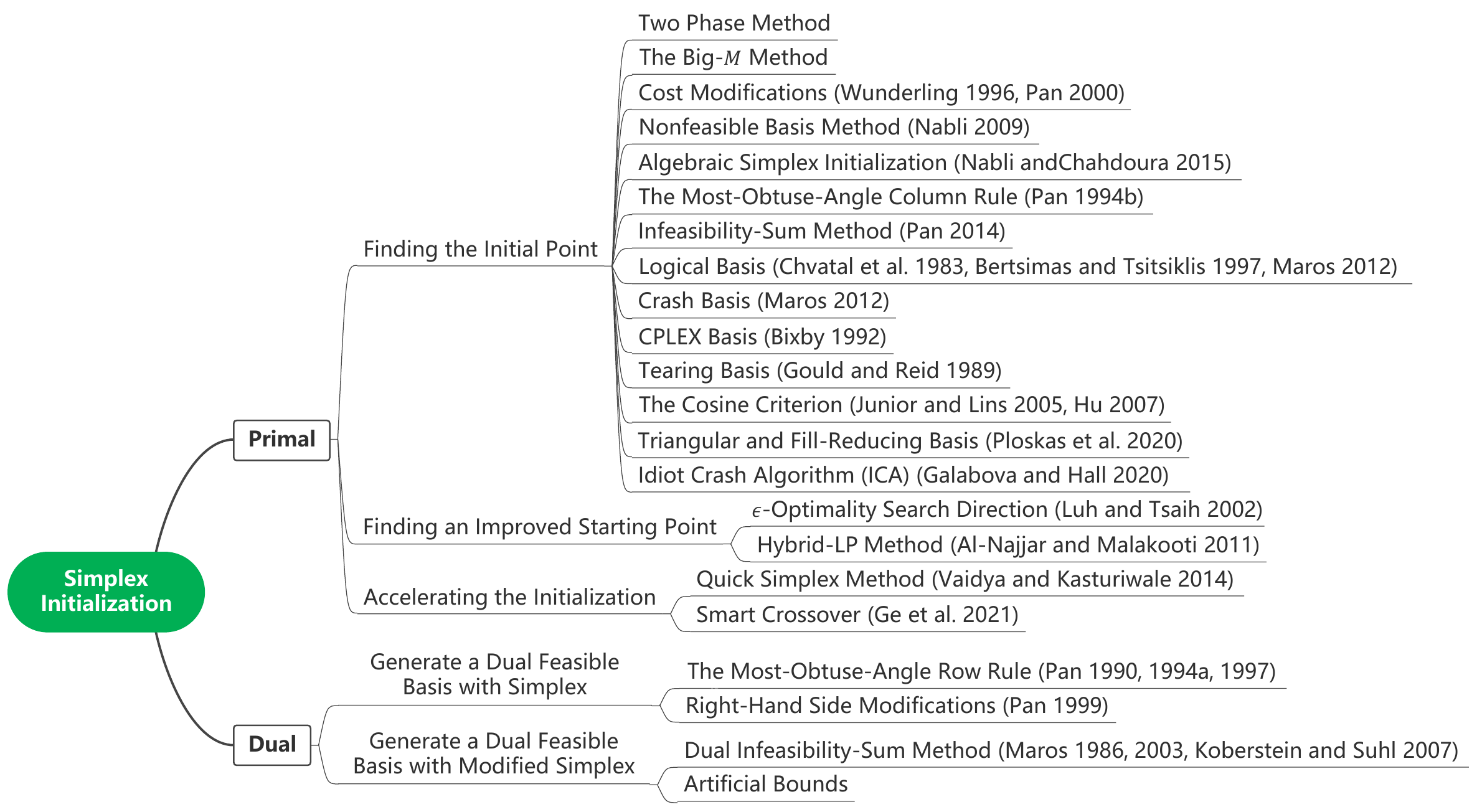}
	\caption{Overview of the initialization methods in simplex}
	\label{fig:overview}
\end{figure}

For the iteration stage, many extension works investigated and proposed new pivot rules to improve the iteration efficiency, e.g., the minimal index pivot rule \citep{bland1977new}, the Edmonds-Fukuda pivot rule  \citep{fukuda1983oriented,clausen1987note}, and the optimal pivot rule \citep{etoa2016new}, etc. For the termination stage, a general survey of different methods is given in~\cite{azlan2020multiplicity}. In our survey, we focus on \textbf{the extension works for improving the initialization stage}.  There are three reasons for studying this stage.  First, the initialization stage plays a significant role in increasing the algorithm efficiency.  A good starting point can lead to less number of iterations or less computation time within each iteration in the following stages, thus resulting in a reduced calculation time of the entire solution process. Second, the computational time to obtain a suitable initial point or basis can sometimes be much greater than the time spent in subsequent stages. Therefore, how to propose some improved algorithms or methods for accelerating the initialization stage is an important research direction.  Third, compared with the iteration stage and the termination stage, there are still many research gaps in the initialization stage that need to be done or completed. An overview of the related simplex initialization methods summarized in this survey is illustrated in Figure~\ref{fig:overview}.

The remainder of the survey is organized as follows. In the Preliminary section, we present the main algorithms of the simplex method and the IPM and compare these two methods. In the Initialization in Primal Simplex section and the Initialization in Dual Simplex section, we investigate the extension works for improving the initialization stage from the perspective of the primal simplex and the dual simplex, respectively. In the final section, we summarize the contributions of this survey and propose some possible future research directions.

\emph{Notations:} For a matrix $A$, $A_{i\bullet}$ and $A_{\bullet j}$ denote the $i$-th row and $j$-th column of $A$, respectively, and $A_{ij}$ represents the element at $i$-th row and $j$-th column in $A$.
${A}^{T}$ and ${A}^{-1}$ respectively denote the transpose and inverse of matrix ${A}$. ${\rm rank}(A)=m$ denotes the rank of matrix $A$.
For a vector $b$, $b_i$ denotes the $i$-th element of $b$. 
$\mathbb{R}$ is the set of real numbers and $\mathbb{R}^n$ is the $n$-dimensional Euclidean space. 
$\mathbb{R}^{m\times n}$ denotes the space of $m\times n$ real matrices. Given two sets $C_1$ and $C_2$, $C_1\backslash C_2=\{s\in C_1| s\notin C_2\}$.
$\cup$ denotes the intersection of sets.
$\parallel\cdot\parallel_2$ and $\mid\cdot\mid$ respectively denote the Euclidean norm of a vector and the absolute value of a scalar. 
${I}_m$ denotes a $m\times m$ identity matrix.

\section{Preliminary}
As mentioned above, the simplex method and the IPM are two main branches of solving LPs. In practice, there are always LPs where one algorithm significantly outperforms the other. In this section, we briefly review these two algorithms and compare them from different perspectives.
\subsection{Linear Programming Formulation}

\subsubsection{Primal/Standard Form.}
Given a general LP problem, it can be formulated into the primal/standard form as:
\begin{equation}
\label{eq:primal}\tag{P}
\begin{aligned}
\min_{x}\quad 	&c^Tx\\
\text{s.t.}\quad 	&Ax=b\\
               		&x\geq0,
\end{aligned}
\end{equation}
where $c\in\mathbb{R}^n$, $b\in\mathbb{R}^m$ and $A\in\mathbb{R}^{m\times n}$ are problem-dependent parameters, and $x\in\mathbb{R}^n$ is the decision variable. Without loss of generality, we assume ${\rm rank}(A)=m$. Though LP problems may appear in other forms, trivial approaches can be applied to transform them into the standard form \cite{dantzig1965linear}. Therefore, it is sufficient to focus on this standard form to better introduce different LP methods.

\subsubsection{Dual Form.}
The associated dual problem of \eqref{eq:primal} is defined as:
\begin{equation}\label{eq:dual}\tag{D}
\begin{aligned}
\max_{y,s}\quad 	&b^Ty\\
s.t.\quad 	&A^Ty+s=c\\
               		&s\geq0,
\end{aligned}
\end{equation}
where $y\in\mathbb{R}^m$ is the (dual) decision variable associated with $x$ and $s\in\mathbb{R}^n$ is the introduced slack variable.

The mathematical relationship between the primal-dual pair is given by the duality theorems as follows.
\begin{thm}[Weak Duality]
Given arbitrary feasible solutions $x$ to \eqref{eq:primal} and $(y,s)$ to \eqref{eq:dual}, we have $c^Tx\geq b^Ty$. 
\end{thm}
%
%\begin{proof}
\noindent
\emph{Proof.} As $x$ and $(y,s)$ are feasible to \eqref{eq:primal} and \eqref{eq:dual}, respectively, we have
\begin{equation*}
\begin{aligned}
Ax=b, 				&\qquad x\geq 0,\\
 A^Ty+s=c, 	&\qquad s\geq 0.
\end{aligned}
\end{equation*}
Thus
\begin{equation*}
c^Tx-b^Ty = ( A^Ty+s)^Tx-(Ax)^Ty=s^Tx \geq 0.
\end{equation*}
%\end{proof}

\begin{thm}[Strong Duality]
If one of the primal-dual pair problems admits an optimal solution, the optimal solution exists for the other problem, and for any optimal solution pair $x^*$ and $(y^*,s^*)$, the duality gap is zero, i.e., $c^Tx^*=b^Ty^*$.
\end{thm}
\subsection{The Simplex Method}\label{chap:simplex_alg}
 The simplex method searches for an optimal solution by visiting adjacent vertices, i.e., basic feasible solutions, of the feasible region. With the elegantly designed entering/leaving basic variables at each iteration, the objective function will monotonically decrease/increase to the optimal value.

\subsubsection{Basic Solutions.}
As ${\rm rank}(A)=m$, $A$ can be permuted into a partitioned-matrix form, i.e., $A = [A_B,A_N]$, where $A_B\in\mathbb{R}^{m\times m}$ is a non-singular sub-matrix of $A$.
\begin{defn}
Any column collection of $A_B$ is called a basis of \eqref{eq:primal}.
\end{defn}

Letting $B$ and $N$ be the associated column indices of $A_B$ and $A_N$, respectively, \eqref{eq:primal} can then be rewritten into the canonical form as follows:
\begin{equation}\label{eq:basis}
\begin{aligned}
\min_{x_B, x_N} \quad 	&c_B^Tx_B+c_N^Tx_N\\
{\rm s.t.\quad} 	&A_Bx_B+A_Nx_N=b\\
               		&x_B,x_N\geq0,
\end{aligned}
\end{equation}
where $c^T=[c^T_B,c^T_N]$ and $x^T=[x^T_B, x^T_N]$ are permuted and partitioned similarly. The basic solution, which satisfies the equality constraints, is obtained by setting non-basic variables to zero. Thus, the primal basic solution based on the current partition is
\begin{equation}\label{eq:basic solution}
\left\{
\begin{aligned}
&x_B=(A_B)^{-1}b\\
&x_N=0
\end{aligned}
\right.,
\end{equation}
Analogously, \eqref{eq:dual} can be written as
\begin{equation}
\begin{aligned}
\max_{y, s_B, s_N} \quad 	    &b^Ty\\
{\rm s.t.\quad} 	&A_B^Ty+s_B=c_B\\
                    &A_N^Ty+s_N=c_N\\
               		&s_B,s_N\geq0.
\end{aligned}
\end{equation}
Since the primal non-basic variables are complementary to the dual basic variables, the associated dual basic solution is obtained by letting $s_B=0$, i.e.,
\begin{equation}\label{eq:dual-basis}
\left\{
\begin{aligned}
&y=(A_B^T)^{-1}c_B\\
&s_B=0\\
&s_N=c_N-A_N^Ty
\end{aligned}
\right.,
\end{equation}
In literature, $\bar{b}\coloneqq A_B^{-1}b$ is called the right-hand side (RHS) coefficient, $\pi\coloneqq (A_B^T)^{-1}c_B$ is the simplex multiplier and $\bar{c}\coloneqq c_N-A_N^T\pi$ is referred to as the reduced cost. Given the basis $A_B$, the basic solution is said to achieve primal feasibility if and only if $x_B\geq 0$, while it achieves dual feasibility if and only if $s_N \geq 0$. Furthermore, if a basic solution is both primal feasible and dual feasible, it is an optimal solution as well. Additionally, a basis is said to be degenerate if there exists any $x_{i\in B}=0$. Degeneracy will cause cycling or stalling in practice, so as to influence the performance of simplex.

\subsubsection{Pivoting.}
Starting with a feasible basis, the simplex method moves from one basis to a neighboring one, i.e., a basis that differs from the previous one by only one element, while preserving the feasibility. The selection of such entering/leaving (basis) variable is called the pivot rule. Geometrically, as the feasible basic solution is associated with a vertex of the feasible region, the simplex method goes through a vertex-to-vertex path to the optimum. After the pivoting operation, the newly generated bases have three features in common:
\begin{enumerate}[1)]
\item Exactly one column of $A_B$ is changed;\label{1}
\item The feasibility is preserved;\label{2}
\item The objective function decreases/increases monotonically.\label{3}
\end{enumerate}

\subsubsection{The Primal Simplex Method.}
According to the type of feasibility preserved during the iteration, the simplex method can be categorized into two classes, i.e., primal simplex and dual simplex.

The primal simplex method is initialized with a primal feasible basis. The feasibility remains within iterations until optimality or unboundedness is detected. Therefore, the primal simplex algorithm can be summarized as Algorithm \ref{alg:primal}. In the algorithm, $e_q$ denotes the unit vector which takes one at position $q$ and zeros otherwise.
% \clearpage
\begin{algorithm}[!htbp]
	\caption{\bf Primal Simplex}
	\label{alg:primal}
	\begin{algorithmic}[1]
		\Require Arbitrary primal feasible basis $A_B$
		\State {\bf Repeat forever}
			\State \hspace{0.5cm} {\bf If} $s_N\geq0$ {\bf then}
				\State \hspace{1.cm}{\bf Break} with optimality
			\State \hspace{0.5cm} {\bf Else}
				\State\hspace{1.cm} Select the entering index $q=\{j\in N:s_j<0\}$\Comment{monotonicity}
				\State \hspace{1.cm} Compute $\Delta x_B = A_B^{-1}A_Ne_q$ 
				% \footnotemark[1]
				and  $t=\left(\max_{i\in B} \frac{\Delta x_i}{x_i}\right)^{-1}$
				\State \hspace{1.cm} {\bf If} $t\leq 0$ {\bf then}
					\State \hspace{1.5cm} {\bf Break} with unboundedness
				\State \hspace{1.cm} {\bf Else}
					\State \hspace{1.5cm} Select the leaving index $p=\arg\max_{i\in B}\frac{\Delta x_i}{x_i}
					$\Comment{feasibility}
				\State \hspace{1.cm} {\bf End If}
			\State \hspace{0.5cm} {\bf End If}
			\State \hspace{0.5cm} Perform pivoting: $B \gets B\cup\{q\}\backslash\{p\}$\Comment{one column is changed}
	\end{algorithmic}
\end{algorithm}

\subsubsection{The Dual Simplex Method.}
Instead of starting with a primal feasible basis, the dual simplex method requires a dual feasible one. Analogously, the dual simplex algorithm can be summarized as Algorithm \ref{alg:dual}.
\begin{algorithm}[!htbp]
	\caption{\bf Dual Simplex}
	\label{alg:dual}
	\begin{algorithmic}[1]
		\Require Arbitrary dual feasible basis $A_B$
		\State {\bf Repeat forever}
			\State \hspace{0.5cm} {\bf If} $x_B\geq0$ {\bf then}
				\State \hspace{1.cm}{\bf Break} with optimality
			\State \hspace{0.5cm} {\bf Else}
				\State\hspace{1.cm} Select the leaving index $p=\{i\in B:x_i<0\}$\Comment{monotonicity}
				\State \hspace{1.cm} Compute $\Delta s_N = (A_B^{-1}A_N)^Te_p$ and $r=\left(\max_{j\in N} -\frac{\Delta s_j}{s_j}\right)^{-1}$
				\State \hspace{1.cm} {\bf If} $r\leq 0$ {\bf then}
					\State \hspace{1.5cm} {\bf Break} with (primal) infeasibility
				\State \hspace{1.cm} {\bf Else}
					\State \hspace{1.5cm} Select the entering index $q=\arg\max_{j\in N} \left(-\frac{\Delta s_j}{s_j}\right)$\Comment{feasibility}
				\State \hspace{1.cm} {\bf End If}
			\State \hspace{0.5cm} {\bf End If}
			\State \hspace{0.5cm} Perform pivoting $B \gets B\cup\{q\}\backslash\{p\}$\Comment{one column is changed}
	\end{algorithmic}
\end{algorithm}

\subsection{Interior Point Method (IPM)}
Different from the simplex methods, IPMs provide an alternative way of solving LP problems. The main idea is that from an initial point in the feasible region, there exists a path across the interior of the polyhedron, along which an optimal point can be reached. That is to say, the LP problem can be solved within a single iteration if the ``right'' direction, i.e., the exact direction from the initial point to the optimal point, is found. Although this is usually impossible, the idea does motivate the emergence of IPMs. 
	
The most remarkable type of IPMs is the path-following IPM, whose idea is to find an optimal solution by following a central path \citep{sonnevend1986analytical}. Considering the primal-dual pair (\ref{eq:primal}) and (\ref{eq:dual}), by duality theory, it is clear that to find optimal solutions for both primal and dual problems, the following system need to be solved:
	\begin{equation}\label{system1}
	\begin{split}
	&Ax = b, x\geq 0,\\
	&A^Ty+s = c, s\geq 0,\\
	&x_j s_j  = 0, j = 1, \ldots, n.
	\end{split}
	\end{equation}
	Here the first line guarantees primal feasibility, the second one guarantees dual feasibility, and the last one, called complementarity condition, is componentwise and guarantees optimality.
	
	By adding logarithmic barrier terms in the objective functions in (\ref{eq:primal}) and (\ref{eq:dual}), we can remove the nonnegative constraints and construct the following logrithmic barrier problem pair:
	\begin{equation}
	\begin{split}
	&\min_{x}\{c^Tx-\mu\sum_{j=1}^{n}\log(x_j): Ax=b\},\\
	&\max_{y, s}\{b^Ty+\mu\sum_{j=1}^{n}\log(s_j): A^Ty+s=c\},
	\end{split}
	\end{equation}
	where $\mu>0$ is called barrier parameter and $\log(\cdot)$ is logarithmic barrier function, which can guarantee the positiveness of the variables $x$ and $s$.
	
	The primal-dual IPM relaxes the complementarity condition and generates a monotone descent sequence $\{\mu_k\}$ with limit $0$. In the $k$-th iteration, $\mu=\mu_k$ is fixed and the approximate solution of the following nonlinear system is obtained by employing Newton's method:
	\begin{equation}\label{system2}
	\begin{split}
	&Ax= b,\\
	&A^Ty+s = c,\\
	&Xs = \mu {\bf 1}_n,
	\end{split}
	\end{equation}
	where $X=\text{diag}(x_1, x_2, \dots, x_n)$ and ${\bf 1}_n=(1,1,\ldots,1)^T\in\mathbb{R}^n$. For each $\mu>0$, the system (\ref{system2}) has a unique solution $\{x(\mu), y(\mu), s(\mu)\}$, and the primal central path and the dual central path are defined as $\{x(\mu): \mu>0\}$ and $\{(y(\mu),s(\mu)): \mu>0\}$, respectively. The optimum will be approximately reached by letting $\mu\rightarrow 0$.

\subsection{Comparison Between Simplex Methods and IPMs}
The competition between the simplex methods and the interior point methods has lasted for a long time. 
	In the following, we compare the simplex methods and (primal-dual) the interior point methods from five aspects.
	
	\begin{enumerate}[1)]
	\item \textbf{Generated solution.} The simplex methods can generate an optimal basic solution, while IPMs generate a sequence of strictly feasible primal and dual solutions and finally produce an $\epsilon$-optimal solution, i.e., a feasible solution pair satisfying $c^Tx-b^Ty<\epsilon$, where $\epsilon$ is small enough. 

	\item \textbf{Complexity \& Practical performance.} 
	Since the number of vertices of the polyhedron might increase exponentially with the problem dimension, and the simplex method might go through all the vertices in the worst case, it has exponential time computational complexity \citep{klee1972good}. Fortunately, various heuristics can be applied to enhance the practical performance of the simplex method. As a result, the simplex method performs much better than its theoretical worst-case bound. Furthermore, the expected complexity of some simplex methods is proved to be polynomial under a probabilistic model \citep{borgwardt2012simplex}. 
	The worst-case complexity of IPMs is polynomial time. Moreover, IPMs are thought to be superior to the simplex method in solving large-scale sparse LP problems by some scholars \citep{andersen1996implementation}.
	
	\item \textbf{Initialization.} 
	In most cases \citep{maros2012computational}, due to the requirement of a feasible basic solution to start with, the simplex methods consist of two phases \citep{dantzig1965linear}, where phase \uppercase\expandafter{\romannumeral1} is used to obtain the required feasible basic solution. A good initial basic solution can greatly reduce the computational cost of Phase II. However, in many cases, Phase I is much more time-consuming than Phase II.
	For IPMs, an initial point, which should be an interior point of the polyhedron, i.e., a feasible point for the LP problem, is needed. Usually, it is difficult to find a feasible initial point in practice. 
	Although infeasible IPM \citep{lustig1990feasibility} can be applied without such a feasible initial point, the price is some weaker complexity. 
	Homogeneous and self-dual method \citep{ye1994nl} might be the most efficient interior-point method so far. It can be initialized by any point whose coordinate components are all positive, and its cost is only a slight increase in the computation of each iteration. 
	
	\item \textbf{Warm-start.} 
	In some applications, after the original LP problem was solved, a new LP problem, which is derived by making some small modifications to the original one (e.g., perturbing the bound of some variables, adding/dropping variables or constraints, etc.), is needed to be solved. By starting from a point (a vertex for the simplex methods and an interior point for IPMs) yielded from the solving process of the original problem, it is expected that fewer steps/iterations are usually required to solve the new modified problem since the obtained start point could be very close to an optimal point. This strategy is called ``warm start''. For the simplex methods, starting a new problem from the optimal point of the original problem allows reaching the optimal point of the new problem fast in most cases. However, despite many efforts that have been made \citep{gondzio1998warm, gondizo1999warm, yildirim2002warm, gondzio2002reoptimization, benson2007exact, gondzio2008new, john2008implementation}, the warm start of IPMs does not show significant efficiency as that of the simplex methods do.
	
	\item \textbf{Generalizations to other optimization problems.} 
	When using cutting plane methods to solve integer LP problems, a basic solution is required to generate cuts \citep{caroe1997cutting, ascheuer1993cutting}. Therefore, the simplex methods have a natural advantage over IPMs. More importantly, due to the efficiency of warm start, dual simplex methods perform extremely well in branch-cut-and-bound schemes where a large number of modified subproblems need to be solved \citep{banciu2011dual, koberstein2008progress, ikura1986computational, koberstein2007progress}.
	Since by some extra purification processes, IPMs can also obtain a basic solution \citep{marsten1989implementation, megiddo1991finding}, cutting plane methods are applicable for IPMs as well \citep{mitchell2000computational, mitchell2003polynomial, elhedhli2004integration, ding2004interior}. However, the warm start procedures of IPMs are not that efficient as the simplex methods, so there is no doubt that the simplex methods are superior to IPMs in solving (mixed) integer LP problems \citep{bixby2012brief}. Besides integer LP problems, some variants of simplex methods are designed for some kinds of nonlinear optimization problems \citep{bazaraa2013nonlinear}, such as linear complementary problems, nonlinear and semi-infinite optimization problems, etc. These kinds of problems can also be solved by applying IPMs \citep{nesterov1994interior}. Moreover, IPMs are more efficient when solving conic-linear optimization problems, especially semi-definite and second-order cone optimization problems \citep{ben1998convex}.
	\end{enumerate}

\section{Initialization in Primal Simplex}
\label{chap:primal_initial}
The initialization methods in different primal simplex algorithms can be classified into three types. The first type is to generate an initial point or basis. The second is to obtain an improved point or basis based on a given point or basis. Then the improved one is utilized as the starting point of the following steps. The third type is to accelerate the calculation process of the first two types. In the following subsections, methods belonging to these three types will be investigated, respectively.

\subsection{Finding the Initial Basis or Point}
\label{find_sub}
\subsubsection{Two-Phase Method.}\label{chap:two-phase}
% \noindent
% \textbf{Two-Phase Method}\\[4pt]
The simplex method usually proceeds in two phases. Phase I is terminated either with a feasible basic solution, or with evidence that the problem is infeasible. If in Phase I, a feasible basic solution is successfully found, then in Phase II, starting from the obtained solution, the simplex algorithm introduced before can be executed in search of the optimum.

Actually, in the two-phase method, Phase I proceeds similarly to Phase II, except that it instead deals with an auxiliary problem, i.e.,
\begin{equation}\label{eq:phase_i}
\begin{aligned}
\min_{x, x_a} \quad 	&{\bf 1}_m^Tx_a\\
{\rm s.t.\quad} 	&Ax + x_a=b,\\
				&x,x_a\geq 0,
\end{aligned}
\end{equation}
where $x_a\in\mathbb{R}^m$ is the introduced artificial variable and ${\bf 1}_m\in\mathbb{R}^m$ denotes the vector of ones. Since for any row with $b_i<0$, we can obtain $b_i>0$ by multiplying both sides by $-1$, without loss of generality, we can assume $b\geq0$. Therefore, the auxiliary problem has a straightforward feasible basic solution, i.e., $x=0$ and $x_a=b$.  Starting from this solution, we can solve \eqref{eq:phase_i} with the simplex algorithm, and encounter two different cases at optimality:

\begin{description}
	\item[Case A]$x_a\neq0$: The original problem \eqref{eq:primal} is infeasible;
   	\item[Case B]$x_a=0$: There are two possibilities:
	\begin{description}
   		\item[\quad Sub-case B.1] No artificial variable remains in the basis: The basis of \eqref{eq:phase_i} is immediately a feasible basis of \eqref{eq:primal};
   		\item[\quad Sub-case B.2]
   		\label{substep}
   		At least one artificial variable remains in the basis: Without loss of generality, assume such variable is in the $i$-th row, where $\bar{b}_i$=0. Select a column $j$ with $(A_B^{-1}A_N)_{ij}\neq0$. Perform the pivoting with $x_j$ as the entering variable and the basic artificial variable in the $i$-th row as the leaving variable. After that, go to Sub-case B.1.
   		
   	% 	Select any $(A_B^{-1}A_N)_{ij}\neq0$ as the pivot and apply pivoting so that the artificial variable is removed from the basis. Therefore, we go to sub-case B.1.
	\end{description}
\end{description}

Though the two-phase method can guarantee a feasible basic solution or evidence of infeasibility at Phase I, it introduces extra artificial variables, thus increasing the dimension, as well as the complexity of the problem. Actually, it has been proved that the problem of determining a feasible solution is of the same complexity degree as solving the LP problem itself \citep{papadimitriou1998combinatorial}. Therefore, Phase I can be very time-consuming in practice, usually even more time-consuming than Phase II \citep{stojkovic2012simplex}.

\subsubsection{Big-$M$ Method.}\label{chap:big-m}
The big-$M$ method is a well-known method to initialize the simplex algorithm. It constructs a feasible basis by introducing artificial variables in the constraints, and eliminates them from the optimal basis by placing large penalty to them in the objective function. Specifically, the auxiliary problem is
\begin{equation}\label{eq:big-m}
\begin{aligned}
\min_{x, x_a} \quad	&c^Tx+M{\bf 1}^T_m x_a\\
{\rm s.t.\quad}	&Ax+x_a=b,\\
					&x,x_a\geq0,
\end{aligned}
\end{equation}
where $x_a\in\mathbb{R}^m$ is the introduced artificial variable and $M\gg0$ is a very large number. Similar to the two-phase method, problem \eqref{eq:big-m} has a trivial feasible basis, i.e., $x_a=b$ and $x=0$. With such a basis, the primal simplex method can be applied to solve the problem. It should be noted that since $M$ is large, high cost will be paid for any $x_a\neq0$. Therefore, though we start with basic variable $x_a=b$, it will be removed from the basis and be pushed to zero in the optimal solution. When $x_a=0$, problem \eqref{eq:big-m} degrades to \eqref{eq:primal} and the obtained solution is directly an optimal solution to \eqref{eq:primal}. If we have $x_a\neq0$ in the solution, the original problem is infeasible.

Actually, the pivoting procedure of big-$M$ is the same as that of the two-phase method. Hence, it is time-consuming as well. However, the introduction of $M$ introduces two more disadvantages. First, it is difficult to determine how large $M$ should be in order to successfully eliminate the artificial variables in practice. Second, numerical issue will occur when we deal with a large number like $M$.

\subsubsection{Cost Modifications.}\label{chap:cost-mod}
The main idea of cost modifications is similar to the two-phase method. It starts with an auxiliary problem which has a straightforward feasible basis, and generates a feasible basis of the original problem by solving the auxiliary one.

Recall the primal problem formulated as \eqref{eq:primal}. The auxiliary problem is constructed by modifying the objective function, i.e.,
\begin{equation}\label{eq:cost-mod}
\hat{c}_{j}=
\left\{
\begin{array}{ll}
A_{\bullet j}^T(A_B^T)^{-1}c_B+\delta_j, 	&\quad\text{if } j\in J\\
c_j,										&\quad\text{otherwise}
\end{array}
\right.,
\end{equation}
where $J=\{j\in N | s_j<0\}$ denotes the set of infeasible variables in the dual problem, and $\delta_j$ is a small positive perturbation to alleviate the problem of degeneracy. It is easy to verify that the basis of \eqref{eq:cost-mod} is dual feasible, i.e., $s_{j\in J}=\delta_j\geq0$. After applying the dual simplex method, we obtain $x_B\geq0$ at optimality. Therefore the optimal basis of \eqref{eq:cost-mod} is immediately a feasible basis to \eqref{eq:primal}, and the primal simplex method can be applied subsequently to compute the optimal solution of the original problem. For a detailed review of this method, as well as its implementation and the corresponding tableau form, readers may refer to \cite{Wunderling1996} and \cite{pan2000primal}.

\subsubsection{Nonfeasible Basis Method.}
% \noindent
% \\ \textbf{Nonfeasible basis method}\\[4pt]
The so-called nonfeasible basis method (NFB) was proposed by Nabli in 2009~\citep{nabli2009overview}. This method is used to construct an initial feasible basis from an infeasible one. 
%Different from other initialization methods based on a nonfeasible basis~\cite{paparrizos2003new,pan2000primal}, 
It can be employed easily without artificial variables and any perturbation in the objective function. 
%or a big $M$ number and without any perturbation in the objective function. 
The feasibility of the obtained basis is achieved by modifying the selection rules of the entering and leaving variables.
This method is completely new and different from the dual simplex algorithm and the criss-cross method~\citep{zionts1969criss}.   
In the same paper~\citep{nabli2009overview}, Nabli introduced the notion of formal tableau. By combining the NFB with the formal tableau, he proposed another new method called formal nonfeasible basis method (FNFB).

As mentioned above, the NFB is used to handle the cases with infeasible initial basis, i.e., the RHS vector $\bar{b}$ has at least one negative component. For such scenarios, it is supposed that the matrix $\beta = A_B^{-1}A_N$ satisfies the following condition:
\begin{equation}
	\forall i\in\{i|\bar{b}_i<0\}, \exists j \text{ s.t. } \beta_{ij}<0.
\end{equation}
If this condition cannot be satisfied, then the LP problem is infeasible. Considering the standard form given in~\eqref{eq:primal}, the procedure of the NFB is given as follows:
%\begin{enumerate}[1.]
%    \item Determine $k=\arg\min\limits_i \{\bar{b}_i | \bar{b}_i<0\}$.
%    \item Build the set $K=\{j | \beta_{kj}<0\}$. \textcolor{red}{If $K=\emptyset$, end the process and the original problem is infeasible.
%    \item Calculate $\omega = -c^T_N + c^T_B A^{-1}_B A_N$.
%   \item Select the entering variable index $q = \arg\min_{j\in K} \{\frac{\omega_j}{\beta_{kj}}\}$.
    % such that $\frac{\omega_q}{\beta_{kq}}\leq\frac{\omega_j}{\beta_{kj}}, \forall j\in K$.
%    \item Select the leaving variable index $p=\arg\max\limits_{i \in \{1,\dots, m\}} \{\frac{\bar{b}_i}{\beta_{iq}} | \bar{b}_i<0 \text{ and } \beta_{iq}<0\}$.
%     \begin{equation}
% 	    p=\arg\max\limits_{i \in \{1,\dots, m\}} \{\frac{\bar{b}_i}{\beta_{iq}} | \bar{b}_i<0 \text{ and } \beta_{iq}<0\}
% 	\end{equation}
%	 \item Repeat the process until $\bar{b}_i \ge 0, \forall i \in \{1,\dots,m\}$ (the feasible basis has been obtained). }

%\end{enumerate}

\begin{enumerate}[1)]
    \item Determine $k=\arg\min\limits_i \{\bar{b}_i | \bar{b}_i<0\}$.
    \item Build the set $K=\{j | \beta_{kj}<0\}$. If $K=\emptyset$, end the process and the original problem is infeasible.
    \item Calculate $\bar{c} = c_N - A_N^T(A_B^T)^{-1}c_B$.
    \item Select the entering variable index $q = \arg\min_{j\in K} \{-\frac{\bar{c}_j}{\beta_{kj}}\}$.
    % such that $\frac{\omega_q}{\beta_{kq}}\leq\frac{\omega_j}{\beta_{kj}}, \forall j\in K$.
    \item Select the leaving variable index $p=\arg\max\limits_{i \in \{1,\dots, m\}} \{\frac{\bar{b}_i}{\beta_{iq}} | \bar{b}_i<0 \text{ and } \beta_{iq}<0\}$.
%     \begin{equation}
% 	    p=\arg\max\limits_{i \in \{1,\dots, m\}} \{\frac{\bar{b}_i}{\beta_{iq}} | \bar{b}_i<0 \text{ and } \beta_{iq}<0\}
% 	\end{equation}
	 \item Repeat the process until $\bar{b}_i \ge 0, \forall i \in \{1,\dots,m\}$ (the feasible basis has been obtained).

\end{enumerate}

% \begin{tabular}{m{12.8cm}}
% 	1. determine $k=\arg\min\{\bar{b}_i | \bar{b}_i<0\}$,\\
% 	2. build the set $K=\{j | \beta_{kj}<0\}$,\\
% 	3. select the entering index $q$ such that $\frac{\pi_q}{\beta_{kq}}\leq\frac{\pi_j}{\beta_{kj}}, \forall j\in K$.\\
% \end{tabular}

% Then, the leaving index $p$ is selected as following:
% \begin{equation}
% 	p=\arg\max_i\{\frac{\bar{b}_i}{\gamma_i} \text{ for } i=1,\ldots,m | \bar{b}_i<0 \text{ and } \gamma_i<0\},
% \end{equation}
% where $\gamma=A_B^{-1}(A_N)_{\bullet q}$.

\subsubsection{Algebraic Simplex Initialization.}
% \noindent
% \\ \textbf{Algebraic simplex initialization} \\[4pt]
In 2015, \cite{nabli2015algebraic} developed a new initialization method based on the notion of linear algebra and Gauss pivoting (hereinafter referred to as algebraic simplex initialization). This method can find a nonsingular initial basis, i.e., $A_B$ is nonsingular, but not necessarily feasible. Therefore, this method was combined with the NFB in the previous subsection by the authors to achieve feasibility. Also, a new pivot rule for the NFB was proposed in this paper~\citep{nabli2015algebraic}, which is advantageous in reducing the number of iterations and the computational time. 

\begin{figure}[htbp]
	\centering
	\includegraphics[width=4in]{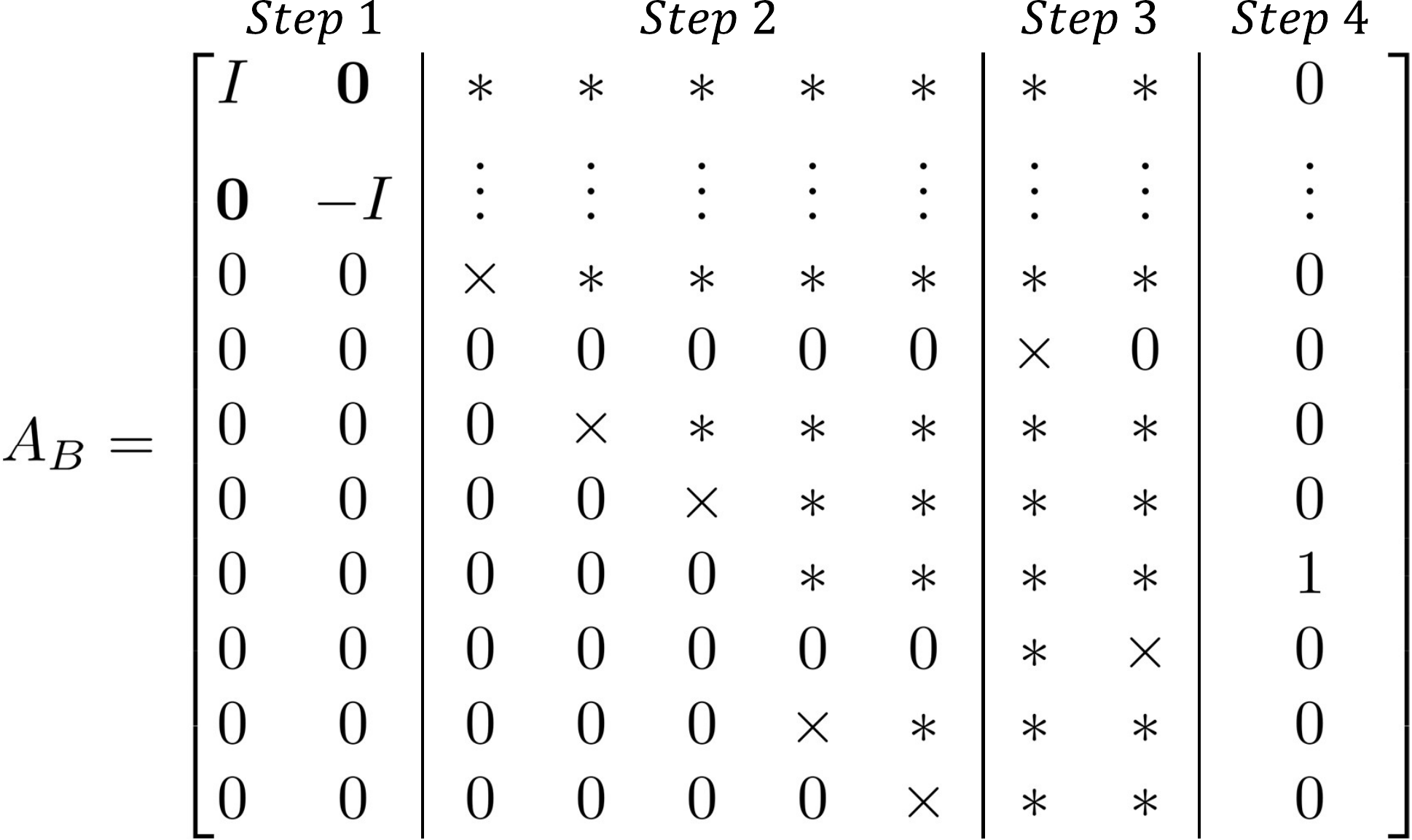}
	\caption{The form of the obtained matrix $A_B$. The symbol `$\times$' indicates that the corresponding entry is nonzero, whereas `$*$' means that the entry is arbitrary.}
	\label{om}
\end{figure}

The algebraic simplex initialization consists of at most four consecutive steps. The first step is to select all the slack variables as basic variables and put their corresponding columns in the coefficient matrix into the formed matrix $A_B$, which is empty before this step. If the LP problem contains no equality constraint, then the obtained basis has been valid, i.e., the obtained $A_B$ is nonsingular; otherwise, the subsequent steps need to be executed. The second and the third step are straightforward. Their main purpose is to continue to select variables from decision variables as new basic variables and fill the columns of the matrix $A_B$ accordingly, so that all the columns of the formed matrix $A_B$ are linearly independent. After these two steps, if the formed matrix $A_B$ is still nonsingular, the last step is required. In order to complete the basis $A_B$, some so-called pivoting variables need to be introduced in this step. Finally, the obtained matrix $A_B$ has the following form as shown in Figure \ref{om}.

Note that the pivoting variables are different from artificial variables since the objective function and the initial problem do not depend on them, and they will be transformed outside the basic variable set by some pivoting steps. Therefore, the algebraic simplex initialization is also artificial-free. Moreover, the redundant constraints and infeasibility can be detected during the pivoting steps.
\subsubsection{The Most-Obtuse-Angle Column Rule.}
The most-obtuse-angle column rule \citep{pan1994variant} combines, to some degree, the work of searching for a feasible basis with the work of searching for an optimal one \citep{wolfe1965composite}. In detail, this method suggests to achieve primal feasibility by iteratively using a modified dual pivot rule. Geometrically, the leaving variable specifies the most obtuse angle with the uphill direction determined by the entering variable. If the uphill direction is close to the direction of the dual objective function, from Figure \ref{fig:observation} we can conclude that basis constructed in this way is more favorable from the perspective of the objective function. The full procedure of this method is shown in the following.
\begin{enumerate}[1)]
\item Select the entering index $p=\arg\min_{i\in B}x_i$. If $x_p\geq0$, the basis is already feasible and go to step \ref{step:final}.\label{step:initial}
\item Compute $\Delta s_N=(A_B^{-1}A_N)^Te_p$. If $\Delta s_N\geq 0$, the algorithm terminates with infeasibility. Otherwise, select the leaving index $p=\arg\min_{i\in B}\Delta s_i$.
\item Perform pivoting $B \gets B\cup\{q\}\backslash\{p\}$ and go to step \ref{step:initial}.
\item Apply primal simplex to compute the optimum. \label{step:final}
\end{enumerate}

Since the feasibility of other variables cannot be maintained in this method, cycling may occur even with no degeneracy. Although \cite{guerrero2005phase} provide a cycling example, such problem may rarely appear in practice.
\subsubsection{Infeasibility-Sum Method.}\label{chap:primal infeasible}
This method is named as infeasibility-sum method because it involves an auxiliary problem which intends to minimize infeasibility-sum, i.e., the negative sum of the infeasible variables. When the infeasible variables are all eliminated in a certain iteration, a feasible basis is obtained.

Let $I=\{i\in B|x_i<0\}$ denote the set of the infeasible basic variables and construct the auxiliary problem as
\begin{equation}\label{eq:primal-infeasible}
\begin{aligned}
\min_{x} \quad		&-\sum_{i \in I}x_i\\
{\rm s.t.\quad}		&A_Bx_B+A_Nx_N=b\\
						&x_{B/I}, x_N\geq0
\end{aligned}
\end{equation}
Note that here we merely impose non-negativity on the variables which already satisfy the non-negative constraints. Therefore, the basic solution \eqref{eq:basic solution} is feasible to \eqref{eq:primal-infeasible} and the primal simplex method can be applied to minimize the infeasibility-sum, so as to compute a feasible basis of \eqref{eq:primal}. It should be noted that as we only impose non-negative constraints on part of the variables, the selection of leaving index is slightly different from the original algorithm.
\begin{thm}
Let $\{\bar{y}_B,\bar{s}_B,\bar{s}_N\}$ be the dual basic solution of the auxiliary problem. The original problem is infeasible if $\bar{s}_N\geq0$.
\end{thm}

Based on the above theorem, the detailed procedure of the infeasibility-sum method is shown below.
\begin{enumerate}[1)]
\item If $x_B\geq0$, the basis is feasible and go to step \ref{final step}. \label{first step}
\item Form an auxiliary problem with respect to the current basis and compute $\bar{s}_N$. If $\bar{s}_N\geq0$, stop with infeasibility. Otherwise, apply one iteration of the modified primal simplex method and go to step \ref{first step}.
\item Apply primal simplex to compute the optimum of the original problem.\label{final step}
\end{enumerate}

For more details about the procedure of this method as well as the modified primal simplex algorithm, readers can refer to~\cite{ping2014linear}.

\subsubsection{Logical Basis.}
% \noindent
% \\ \textbf{Logical basis}\\[4pt]
Logical basis is the simplest initial basis~\citep{chvatal1983linear, bertsimas1997introduction, maros2012computational}. To form such a basis, all constraints (equality \& inequality) add a distinct logical variable after the decision variables. Then, the corresponding columns vectors of all the logical variables form a unit matrix that can be used as an initial basis, i.e., $A_B=I$.

The logical basis has three main advantages: first, its creation is trivial; second, the inverse of $A_B$ is just identity matrix $I$, which is available without any computation; third, the first iterations are very fast as the LU factorization is sparse. 
However, the logical basis generally leads to substantially more iterations, and thus more advanced initial bases are expected.

\subsubsection{Crash Basis.}
% \noindent
% \\ \textbf{Crash basis}\\[4pt]
Compared with an extremely sophisticated algorithm which can provide a good initial basis, a crude algorithm which may provide a reasonably good initial basis quickly is more favored. 
Therefore, some heuristic algorithms, called crashing, have emerged. 
These crash algorithms are used to quickly find a good initial basis with as many decision variables as possible. Usually, the obtained basis is a triangular basis due to some irreplaceable benefits: first, there will be no fill-in in the subsequent iterations; second, it is numerically accurate to calculate the inverse of a triangular matrix; third, it is easy to create a triangular basis; last but not least, operations with triangular matrices are less time-consuming. 
In the LP context, there are two types of triangular basis: the lower triangular basis and the upper triangular basis. Both of the two types have zero-free diagonals.

%The conceptual framework of most triangular crash procedures is as follows.
Most triangular crash procedures are based on the same conceptual framework as following. 
First, partition the coefficient matrix $A$ into $[\hat{A},I]$, where $\hat{A}$ corresponds to the decision variables and $I$ corresponds to the logical variables. Then, define the row and column counts, $R_i$ and $C_j$, as the number of nonzeros in the $i$-th row of $\hat{A}$ and that in the $j$-th column of $\hat{A}$, respectively. 
For lower triangular basis, select a pivot row $i=\arg\min_k\{R_k\}$. If $R_i=1$, the pivot column is unique; otherwise, the column with the smallest column count should be selected to keep the number of nonzero elements in $A_B$ small. To avoid the transformation at inversion/factorization, all the columns with nonzero in row $i$ will not be considered in the subsequent selection process. Then, update row and column counts for the remaining rows and columns of $\hat{A}$ and repeat the above procedure. The main idea for the upper triangular basis is similar. 

In practice, the triangular crash procedures include more other considerations, such as feasibility (CRASH(LTSF)) and degeneracy (CRASH(ADG)). Some of them can be found in~\cite{maros2012computational, maros1998strategies} for more details.

\subsubsection{CPLEX Basis.}
% \noindent
% \\ \textbf{CPLEX basis}\\[4pt]
CPLEX basis was proposed by~\cite{bixby1992implementing}. The essential purpose is to construct a sparse, well-behaved basis with as few artificial variables and as many free variables as possible. Therefore, the CPLEX basis may not include many variables in any optimal basis, but can reduce the work of removing artificial variables. 

To find such a CPLEX basis, a preference order of the variables should be constructed first. Consider the given LP problem in the following form:
\begin{equation}
	\begin{split}
		\min_{x, s_1, s_2} \quad &c^Tx\\
		\text{s.t.}\quad 
		& A_1 x + s_1 = b_1,\\
		& A_2 x - s_2 = b_2,\\
		& A_3 x = b_3,\\
		& l\leq x\leq u,\\
		& s_1\geq 0,\quad s_2\geq 0,
	\end{split}
\end{equation}
where $x=(x_1, \ldots, x_n)^T$ are decision variables, and $s_1 = (x_{n+1}, \ldots, x_{n+m_1})^T$ and $s_2 = (x_{n+m_1}, \ldots, x_{n+m_1+m_2})^T$ are slack variables. All the indices of variables can be divided into four sets:
\begin{equation*}
	\begin{split}
		&C_1=\{n+1, \ldots, n+m_1+m_2 \},\\
		&C_2=\{j:x_j \text{ free} \},\\
		&C_3=\{j\leq n: \text{ exactly one of } l_j \text{ and } u_j \text{ is finite} \},\\
		&C_4=\{j:-\infty\leq l_j,u_j\leq +\infty\},
	\end{split}	
\end{equation*}
where $C_i$ will be prefered to $C_{i+1}$  $(i=1,2,3)$. Note that $C_1$ is just the set of indices of all the slack variables, which is the most prefered set due to the sparsity and numerical properties. 

Then, define a penalty $\bar{q}_j$ for $j\in\{1,\ldots,n+m_1+m_2\}$:
\begin{equation}
	\bar{q}_j=\left\{
	\begin{aligned}
		&0, \quad\text{if } j\in C_2,\\ 
		&l_j, \quad\text{if } j\in C_3 \text{ and } u_j=+\infty,\\
		&-u_j, \quad\text{if } j\in C_3 \text{ and } l_j=-\infty,\\
		&l_j-u_j, \quad\text{if } j\in C_4.
	\end{aligned}
	\right.
\end{equation}
Let $c=\max\{|c_j|:1\leq j \leq n\}$ and define
\begin{equation}
	c_\text{max}=\left\{
	\begin{aligned}
		&1000c, \quad\text{if } c\neq 0,\\ 
		&1, \quad\text{otherwise}.
	\end{aligned}
	\right.
\end{equation}
Finally, for $j\in\{1,\ldots,n\}$, define
\begin{equation}
	q_j=\bar{q}_j+c_j/c_\text{max}.
\end{equation}
The indices in sets $C_2$, $C_3$ and $C_4$ are sorted in ascending order of $q_j$. The lists are concatentated into a single ordered set $C=\{j_1,\ldots,j_n\}$ with the freest decision variable in the front. 
Now, the basis $A_B$ can be constructed according to the steps as shown in~\cite{bixby1992implementing}.

The construction of the CPLEX basis is quite simple and fast. As a result, it is considered a good choice of the default initial basis. The computational results suggested that the CPLEX basis has good performance on easy problems, but it is generally less effective for difficult ones.

\subsubsection{Tearing Algorithm.}
% \noindent
% \\ \textbf{Tearing algorithm}\\[4pt]
%Gould and Reid~
\cite{gould1989new} proposed a remarkable initialization algorithm for large-scale and sparse LP problems, which can find an initial basis as feasible as possible with a reasonable computational cost. This algorithm is called tearing algorithm. Its main idea is to break the initialization problem into several smaller pieces and solve each of them. 
There are two main assumptions for the tearing algorithm. First, the coefficient matrix $A$ can be transformed into a lower block triangular matrix with small blocks by permutating its rows and columns. Second, an efficient simplex solver is available for solving dense LP problems with fewer than $t$ rows, where $t$ is a small number. It is assumed that after some row and column permutations, the permuted matrix has the following form:
\begin{equation}\label{pa}
\begin{bmatrix}
  A_{11}& & & & \\
  A_{21}&A_{22}& & & \\
  \vdots&\vdots&\ddots& & \\
  A_{r1}&A_{r2}&\cdots&A_{rr}& \\
  A_{(r+1)1}&A_{(r+1)2}&\cdots&A_{(r+1)r}&0\\
 \end{bmatrix}
\end{equation}
where $A_{ij}\in\mathbb{R}^{m_i\times n_j}$. Note that $m_i$ and $n_j$ are positive integers for all $i,j\in\{1,2,\ldots,r\}$, but $m_{r+1}$ may be zero and thus nonnegative. Generally, the size $m_i$ of the blocks are very small. Such a block lower-triangular form can be obtained by the algorithm proposed by~\cite{erisman1985structurally}. Then, the following problem is considered:
\begin{subequations}
	\begin{align}
		\min_{v,w} \quad&{\bf 1}_m^T(v+w) \label{t1a}\\
		\text{s.t.}\quad &
		\begin{bmatrix}
A_{11} & & & & \\
A_{21} & A_{22} & & & \\
\vdots & \vdots & \ddots & & \\
A_{r 1} & A_{r 2} & \cdots & A_{r r} & \\
A_{(r+1) 1} & A_{(r+1) 2} & \cdots & A_{(r+1) r} & 0
\end{bmatrix}\begin{bmatrix}
x_{1} \\
x_{2} \\
\vdots \\
x_{r} \\
x_{r+1}
\end{bmatrix}+\begin{bmatrix}
v_{1} \\
v_{2} \\
\vdots \\
v_{r} \\
v_{r+1}
\end{bmatrix}-\begin{bmatrix}
w_{1} \\
w_{2} \\
\vdots \\
w_{r} \\
w_{r+1}
\end{bmatrix}=\begin{bmatrix}
b_{1} \\
b_{2} \\
\vdots \\
b_{r} \\
b_{r+1}
\end{bmatrix},\label{t1b}\\
		& l_i\leq x_i \leq u_i, \quad v_i\geq 0, \quad w_i\geq 0, \quad i=1,\ldots,r+1,\label{t1c}
	\end{align}
\end{subequations}
where $x_i \in \mathbb{R}^{n_i}$, $v_i,w_i,b_i \in\mathbb{R}^{m_i}$ and $l_i$ and $u_i$ are the lower bound and upper bound of $x_i$, respectively. 
Note that a feasible solution can be reached with ${\bf 1}_m^T(v+w)=0$. Although sometimes the feasibility may not be achieved, a basis near to the feasible region can be obtained. 

It is easy to know that the first block of \eqref{t1b} requires that the following conditions are satisfied:
\begin{equation}\label{c1}
	A_{11}x_1+v_1-w_1=b_1, \quad l_1\leq x_1 \leq u_1,\quad v_1,w_1\geq 0,
\end{equation}
and it is expected that both $v_1$ and $w_1$ are driven to zero.
Assuming that $m_i\leq t, \forall i$, this may be achieved by using DLP to minimze ${\bf 1}^T_{m_1}(v_1+w_1)$ subject to \eqref{c1}, which produces the optimal solution $\hat{x}_1, \hat{v}_1, \hat{w}_1$ and a set of $m_1$ basic variables probably including some of $v_1$ and $w_1$. The obtained solution $\hat{v}_1$ and $\hat{w}_1$ will be zero if the original problem is feasible.

Moving to the $k$-th stage ($1<k\leq r$) of this algorithm, the optimal solution $\hat{x}_i, \hat{v}_i, \hat{w}_i, \forall 1\leq i <k$ and a set of $m_1+\cdots+m_{k-1}$ basic variables have been obtained. Then, the following conditions need to be satisfied:
\begin{equation}\label{ck}
	A_{kk}x_k+v_k-w_k=b_k-\sum_{i=1}^{k-1}A_{ik}\hat{x}_i, \quad l_k\leq x_k \leq u_k,\quad v_k,w_k\geq 0,
\end{equation}
and $v_k$ and $w_k$ had better to be zero. Similarly, this can be achieved by minimizing ${\bf 1}^T_{m_k}(v_k+w_k)$ subject to \eqref{ck}. 
Last, as the $(r+1)$-th block is $0$, $x_{r+1}$ will have no contribution. Therefore, $\hat{v}_{r+1}$ and $\hat{w}_{r+1}$ should be calculated according to
\begin{align}
	&\hat{v}_{r+1}=\max\{0,b_{r+1}-\sum_{i=1}^{r}A_{(r+1)i}\hat{x}_i\},\\
	&\hat{w}_{r+1}=\max\{0,-b_{r+1}+\sum_{i=1}^{r}A_{(r+1)i}\hat{x}_i\},
\end{align}
where the maximum is taken by elements. Variables in $\hat{v}_{r+1}$ and $\hat{w}_{r+1}$ with nonzero value are chosen to be basic. To make up a full basis, some variables should be arbitrarily picked to cover the remaining rows at the end of the process. 

For all $k>1$, if $\hat{v}_k\neq 0$ or $\hat{w}_k\neq 0$, the backtracking will be executed: if $m_j+m_{j+1}+\cdots+m_k\leq t$ holds for some $j$, then the subproblems in stage $j, j+1, \ldots, k$ can be integrated as one problem, which can be solved by DLP. 

\subsubsection{The Cosine Criterion.}
\begin{figure}[!htbp]
		\centering
		\includegraphics[]{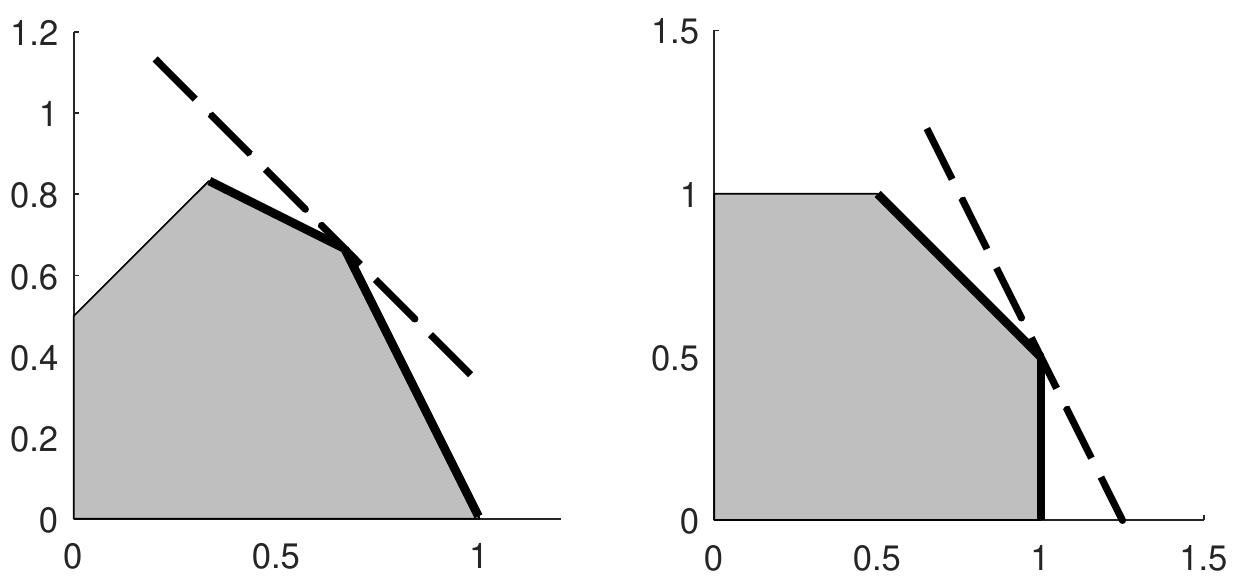}
		\caption{Illustration of the observation, plotted using Plot 2D/3D region \citep{plotrange2021}. The bold lines represent the constraints that make the minimum angle with the objective function, which is denoted by the dashed line.}
		\label{fig:observation}
	\end{figure}

% \noindent 
% \\ \textbf{\textcolor{red}{The Cosine Criterion}}\\[4pt]
The cosine criterion is inspired by the observation that the optimal vertex is usually formed by the constraints that make the minimum angle with the objective function (Figure \ref{fig:observation}). Though similar idea has been researched in~\cite{cos1993,stojkovic2001two}, these algorithms cannot be implemented efficiently due to the existence of redundant constraints.~\cite{junior2005improved} and~\cite{hu2007note} proposed new algorithms which can handle the redundant constraints. In these algorithms, though the cosine criterion cannot guarantee an optimal solution, the obtained vertex turns out to be a near-optimal point. Starting from such a vertex can reduce the number of iterations required by the simplex method, thus speeding up the solution process.

% \footnotetext{Per Bergstr$\ddot{\rm o}$m (2021). \href{https://www.mathworks.com/matlabcentral/fileexchange/9261-plot-2d-3d-region}{Plot 2D/3D region}, MATLAB Central File Exchange. Retrieved April 28, 2021.}

With a bit of abuse of notations, we initialize $B=\varnothing$ and let $N=\{1,\dots,n\}$ be the corresponding complementary set. At each time, one variable is moved from $N$ to $B$, i.e., $B=B\cup\{q\}$ and $N=N\backslash\{q\}$, where $q$ is selected based on the angle and the rank of $A_B$, i.e.,
\begin{equation}
%k_q = \arg\max_j\{\alpha_j| \left\| \bar{N_2} \right\|_{\infty}\neq0\},
q = \arg\max_{j\in N}\{\alpha_j| \left\| \bar{(N_2})_{\bullet j} \right\|_{\infty}\neq0\},
\end{equation}
%where $\alpha_j=(A_{\bullet j})^Tb\big/\left\|A_{\bullet j}\right\|$ is in proportion to the cosine of the angle between $A_{\bullet j}$ and $b$, and $\bar{N}_2$ is a matrix calculated based on LU factorization. The condition $\left\| \bar{N_2} \right\|_{\infty}\neq0$ ensures that the constructed basis $A_B$ is non-singular. According to the feasibility of the obtained basis, either primal simplex or dual simplex is applied to solve the problem. Nevertheless, if the basis is infeasible, other initialization methods introduced in this section can be performed to generate a feasible one.
where $\alpha_j=(A_{\bullet j})^Tb\big/\left\|A_{\bullet j}\right\|$ is named as the dual pivoting index, which is proportional to the cosine of the angle between $A_{\bullet j}$ (the $i$-th constraint of \eqref{eq:dual}) and $b$ (the objective function of \eqref{eq:dual}), and $\bar{N}_2$ is a matrix calculated based on LU factorization. The condition $\left\| \bar{N_2} \right\|_{\infty}\neq0$ ensures that the constructed basis $A_B$ is non-singular. According to the feasibility of the obtained basis, either primal simplex or dual simplex is applied to solve the problem. Nevertheless, if the basis is infeasible, other initialization methods introduced in this section can be performed to generate a feasible one.

The advantage of the cosine criterion is that it can reduce the number of iterations significantly, up to 40\% on Netlib problems. However, the calculation of $\bar{N_2}$ requires LU factorization, which unfortunately tends to be time-consuming as well. Furthermore, as the obtained basis is not likely to be sparse, the computation time per iteration may increase. Therefore, the overall efficiency may not be improved much.

\subsubsection{Triangular and Fill-Reducing Basis.}
% \noindent
% \\ \textbf{Triangular and Fill-reducing Basis} \\[4pt]
One computation difficulty of the simplex method lies in the calculation of basis inverse. In computational practice, LU factorization is used to accelerate such process. In order to further improve the time efficiency, %Plosksa et al.
\cite{ploskas2020triangulation} intended to find a sparse and near-triangular basis so that the factorization becomes easier. In this case, though the number of iterations may increase, the computation time per iteration is reduced and the overall efficiency is improved.

Ignoring the object and the constraints, the algorithm permutes matrix $A$ as:
\begin{equation}
A = \left[
	\begin{matrix}
	A_{11}	& A_{12}\\
	0       	& A_{22}
	\end{matrix}
\right],
\end{equation}
where $A_{11}$ represents the maximal diagonal factors in $A$. If ${\rm rank}(A_{11})\geq m$, i.e., $A_{11}$ is large enough to form an initial basis, the algorithm stops. Otherwise, $A_{22}$ is subsequently ordered by a fill-reducing order, and the first $m-{\rm rank}(A_{11})$ columns are selected in completion of the basis.

Note that with $A_{11}$, the constructed basis will be as triangular as possible. Additionally, as $A_{22}$ is ordered based on fill-in effect, the sparsity will be preserved during iterations. Therefore, computation time per iteration is highly reduced. However, since both the objective function and the constraints are completely ignored, the initial vertex may be far from optimal, thus requiring more iterations to terminate.

Compared with the cosine criterion, which focuses on a near-optimal starting point, there exists a trade-off between fast iteration and less iteration times. Near-optimal basis tends to require less iterations, while sparse and near-triangular basis speeds up the iteration process itself.

\subsubsection{Idiot Crash Algorithm (ICA).}
\label{ICA_subsection}
% \noindent 
% \\ \textbf{Idiot Crash Algorithm (ICA)}\\[4pt]
The main idea of ICA~\citep{galabova2020idiot} is to relax the original LP problem to an approximate problem with ``soft''  constraint, and then solve this relaxed problem to obtain a near-optimal point. This point is later used as the starting point of the simplex method for solving the original problem. 

Recall the standard LP defined in \eqref{eq:primal}, ICA obtains the relaxed problem by replacing the equality constraint with two extra terms in the objective function, i.e., a linear Lagrangian term and a quadratic penalty
term, as follows:
\begin{equation}
\begin{aligned}
\min_{x} \quad 	&c^Tx + \lambda^T(Ax-b) + \frac{1}{2\mu}\|A x-{b}\|_{2}^{2} \\
{\rm s.t.\quad} 	&x \ge 0,
\end{aligned}
\end{equation}
where $\lambda$ is the Lagrange multiplier and $\mu$ is a penalty weight. This relaxed problem can be easily solved by existing methods, such as IPMs. As $\lambda$ goes to infinite, the optimal solution of this relaxed problem will converge to the optimal solution of the original LP problem. To obtain a near-optimal point of the original problem, in each iteration, ICA updates the parameters $\lambda$ and $\mu$ by some heuristic rules and then solves the corresponding relaxed problem. The total iteration number of ICA is finite and predefined heuristically.

As shown in the previous part, the basic simplex method begins with a basic feasible solution (a vertex of the feasible region). However, after finite iterations, the near-optimal point obtained by ICA may be an interior point of the feasible region.
%or an infeasible point
To obtain a basic feasible solution near the point given by ICA, a crossover procedure is added. For some LP problems with special structures, the crossover procedure can be further accelerated to improve the efficiency. These methods will be discussed in the Accelerating the Initialization subsection.

The advantage of ICA is that it can transforms the original problem into a more tractable problem without the equality constraint and obtain a near-optimal point quickly. Nevertheless, one challenge of ICA is how to design the heuristic parameter updating rules. If the rules are not chosen properly, the obtained point will not be a good starting point and the efficiency of the algorithm will not be improved much.

\subsection{Finding an Improved Starting Point}
\label{improve_sub}
Instead of finding an initial basis or point as discussed in the previous subsection, some methods tend to utilize the idea of IPMs to obtain an improved point based on a given initial point, and then regard this new point as the starting point of the simplex method.

\subsubsection{$\epsilon$-Optimality Search Direction.}
% \noindent 
% \\ \textbf{$\epsilon$-Optimality Search Direction}\\[4pt]
The $\epsilon$-optimality search direction algorithm was proposed in~\cite{luh2002efficient}. This algorithm was motivated by the fact that the IPM can approach the neighborhood of the optimal solution faster than the simplex method. In this algorithm, an improved point is obtained by moving in a proposed direction. This point is later used as the starting point of the simplex method for calculating the optimal solution of a given LP problem. The proposed direction combines an interior direction of the feasible region and the negative direction of the objective. 

This algorithm focuses on a normalized LP problem as follows:
\begin{equation}\label{eq:normal_lp}
\begin{aligned}
\min_{x} \quad 	&c^T x\\
{\rm s.t.\quad} 	&Ax \ge b\\
              		&x\geq0,
\end{aligned}
\end{equation}
where $\|c\|^2_2 = 1, \| A_{i\bullet} \|^2_2 = 1, \forall i \in \{1,2,\dots,m\}$. Any LP problem can be easily transferred into the normalized version by choosing $c = \frac{c}{\| c\|^2_2}, A_{i\bullet} = \frac{A_{i\bullet}}{\| A_{i \bullet} \|^2_2}, b_i = \frac{b_i}{\| A_{i\bullet} \|^2_2}$. This normalization process will not change the optimal solution of the original LP problem.

\begin{defn}
Given a feasible point $x$, if $\forall\,\delta >0$, the set $\{x'\mid \|x'-x\|^2_2 < \delta\}$ does not belong to the feasible region, then $x$ is called a boundary point.
\end{defn}

Given a boundary point $x$, this algorithm defines two sets as follows:
\begin{equation}
\begin{aligned}
\Omega_1 &= \{i \mid A_{i\bullet} x =b_i\},\\
\Omega_2 &= \{ j \mid x_j=0\}.
\end{aligned}
\end{equation}
These two sets collect the indices of active constraints at the given boundary point $x$. Based on these two sets, the algorithm calculates a vector $h$ as follows:
\begin{equation}
\label{h_operation}
h = \frac{\sum_{i \in \Omega_{1}} {A}^T_{i\bullet}+\sum_{j \in \Omega_{2}} {e}_{j}}{\left\|\sum_{i \in \Omega_{1}} {A}^T_{i\bullet}+\sum_{j \in \Omega_{2}} {e}_{j}\right\|_2},
\end{equation}
where $e_j$ is a unit vector whose $j$-th element is $1$ and other elements are $0$. The dimension of $e_j$ is consistent with $A^T_i$. Then based on $h$, the direction for obtaining an improved point starting from the given point $x$ is defined as:
\begin{equation}
\label{search_direction}
{g}=\left\{\begin{array}{cl}
0 & \text { if } h=c, \\
Proj\left(\frac{{h}-{c}}{\left\|h-c\right\|_2}\right) & \text { if } h \neq c,
\end{array}\right.
\end{equation}
where $c$ is the vector in the objective function. The projection operation $Proj(\cdot)$ is used to guarantee the feasibility of the direction and its mathematical form can be found in~\cite{luh2002efficient}. Starting from the current boundary point $x$, the direction $g$ actually points to the interior of the feasible region.

Since the LP problem attempts to minimize the objective, this algorithm also constructs a proper step size $\eta$ and proves that the objective can be reduced when moving in the direction of $g$ with the step size $\eta$. Besides, with this specific step size, the new point obtained after one iteration, i.e., $x' = x + \eta g$, is also a boundary point of the feasible range. This iteration can then be repeated based on this new point $x'$.

One important detail in implementing this algorithm is that the search direction given in~\eqref{search_direction} will be replaced by $Proj(c)$ when the step size is less than a predefined value. This is to improve the efficiency of the algorithm when the current point is close to the optimal point. According to the experimental results given in~\cite{luh2002efficient}, the $\epsilon$-optimality search direction algorithm can reduce the iteration number of the basic simplex method by about $40\%$. 

An extension work of this algorithm was introduced in~\cite{chaderjian2003comments}. In this extension work, a special example is given to show that the denominator term in~\eqref{h_operation} can be zero. Therefore, in order to handle the anomaly where the denominator is 0, the new algorithm changes the definition of $h$ as follows:
\begin{equation}
\label{new_h_operation}
h = 
\begin{cases}
0, & \sum_{i \in \Omega_{1}} {A}^T_{i\bullet}+\sum_{j \in \Omega_{2}} {e}_{j}=0, \\
\frac{\sum_{i \in \Omega_{1}} {A}^T_{i\bullet}+\sum_{j \in \Omega_{2}} {e}_{j}}{\left\|\sum_{i \in \Omega_{1}} {A}^T_{i\bullet}+\sum_{j \in \Omega_{2}} {e}_{j}\right\|_2}, &  \sum_{i \in \Omega_{1}} {A}^T_{i\bullet}+\sum_{j \in \Omega_{2}} {e}_{j} \neq 0.
\end{cases}
\end{equation}

The new algorithm also shows that if the initial direction is chosen as $-c$, the algorithm efficiency can be further improved. In addition, the extension work corrects several mathematical errors in~\cite{luh2002efficient}. 

The $\epsilon$-optimality search direction algorithm and its improved version given in~\cite{luh2002efficient, chaderjian2003comments} can be regarded as auxiliary tasks before the basic simplex method. An interesting issue to be investigated is the trade-off between the number of iteration steps to find an improved starting point and the total computation time to solve the LP problem.

\subsubsection{Hybrid-LP Method.}
% \noindent 
% \textbf{Hybrid-LP}\\[4pt]
Hybrid-LP method was introduced in~\cite{al2011hybrid}. The idea of the hybrid-LP is similar to the $\epsilon$-optimality search direction algorithm, but it differs in that instead of designing the iterative direction to acquire the improved starting point based on a given point, it obtains the direction according to the non-basic variables (NBVs). The hybrid-LP method experimentally shows a reduction in both the iteration number and the total computation time. 

The process of the hybrid-LP method can be divided into the following five steps:
\begin{enumerate}[1)]
	\item Select $k$ NBVs to construct the iteractive direction. The value of $k$ is chosen as:
	$$
	k = \alpha \min(m, n-m),
	$$
	
	where $m$ is the number of constraints and $n-m$ is the number of NBVs. The variable $\alpha \in [0,1]$ is a predfined parameter. 
	
	\item Divide these $k$ selected variables into two sets based on whether a change in the variable will result in an increase or decrease in the objective function. Denote these two sets by $s_I$ and $s_D$, respectively.
	
	\item Construct the iterative direction as $d = \zeta d_I + (1-\zeta) d_D$ where $d_I$ is generated by NBVs in $s_I$, while $d_D$ is given by NBVs in $s_D$. The parameter $\zeta$ is selected based on the rule that the direction can lead to an improved objective value.
	
	\item Given the current point $x$, find the maximum step size $\theta$ such that the point $x'$ obtained after one iteration along the direction $d$, i.e., $x' = x + \theta d$, is still within the feasible range.
	
	\item Find a nearby basic feasible solution (vertex) of $x'$ by the reduction process or some crossover methods, and treat it as the starting point of the basic simplex method to solve the original LP problem.
\end{enumerate}

In the experiments, the parameters $\alpha, \zeta$, etc, are chosen heuristically. Therefore, one possible way to improve the hybrid-LP method is to do an optimization on the parameter selection. Similar to the $\epsilon$-optimality search direction algorithm, it is also meaningful to study the trade-off between the iteration number of the hybrid-LP method and the total running time.

\subsection{Accelerating the Initialization}
\label{acc_subsection}
The methods in the previous two subsections either attempt to obtain an initial point (basis) or try to obtain an improved point along a defined direction. In fact, there are some other methods that try to speed up the efficiency of some of the methods mentioned above, i.e., accelerate the Phase I to improve the efficiency of the overall algorithm for solving the LP problem.

\subsubsection{Quick Simplex Method.}
The idea of the quick simplex method~\citep{vaidya2014quick} is to perform the pivoting based on multiple pairs of variables instead of only one pair. One important application of this method to accelerate the simplex initialization can be shown in the two-phase method introduced in the Finding the Initial Basis or Point subsection~\citep{vaidya2016application}. With the quick simplex method, the pivoting operation of the Phase I in the two-phase method is modified as follows:

\begin{enumerate}[1)]
	\item Define the largest number of variable pairs selected in each pivoting operation as $N_p$.
	\item Calculate the simplex tableau with the current basis.
	\item Select at most $N_p$ entering variables based on the reduced cost given in the simplex tableau and calculate the ratios.
	\item Select the corresponding number of leaving variables with the smallest ratios. The selected leaving variables should not be in the same row of the simplex tableau.
	\item End the process if the remaining basis do not have any artificial variables.
\end{enumerate}

The quick simplex method can also be implemented to accelerate the pivoting of Phase II in the basic simplex method, or other simplex initialization methods with a similar iterative process.
 
\subsubsection{Smart Crossover.}
Most of simplex methods begin with a basic feasible point (vertex) to solve the LP problem. To obtain such a starting point, the crossover operation is needed in some simplex initialization methods, such as the ICA mentioned in the Idiot Crash Algorithm (ICA) subsection. However, the crossover operation can be very time-consuming. Therefore, it is necessary and meaningful to propose some smart crossover methods which are more efficient.

%Ge et al.~
\cite{ge2021interior} introduced two network crossover methods which can deal with the minimum cost flow (MCF) problem and the optimal transport (OT) problem. These two problems are two special types of LPs. Specifically, the MCF problem can be easily transformed into an equivalent OT problem.

\emph{Column Generation Method}:
Given a directed graph $\mathcal{G}=(\mathcal{N},\mathcal{A})$, $\mathcal{N}$ and $\mathcal{A}$ are the entire node set and the arc set, respectively. The node set $\mathcal{I}(i), \forall i\in \mathcal{N}$ includes all nodes that have an arc from the node $i$, while the node set $\mathcal{O}(i), \forall i\in \mathcal{N}$ is comprised of all nodes that have an arc pointing to the node $i$. The general form of the MCF problem is given as follows:
\begin{equation}
\label{MCF_problem}
\begin{aligned}
\min_{x} \quad 	&\sum\limits_{(i, j) \in \mathcal{A}} c_{i j} x_{i j} \\
{\rm s.t.\quad} 	&b_{i}+\sum\limits_{j \in \mathcal{I}(i)} x_{j i}=\sum\limits_{j \in \mathcal{O}(i)} x_{i j}, \forall i \in \mathcal{N} \\
& 0 \leq x_{i j} \leq u_{i j}, \forall(i, j) \in \mathcal{A},
\end{aligned}
\end{equation}
where $c_{i,j}$ and $u_{i j}$ denote the cost and the capacity limit on the arc $(i, j)\in \mathcal{A}$, respectively. The variable $b_i$ represents the external supply of the node $i \in \mathcal{N}$. 

The goal of the MCF problem is to design the amount of flow on each arc, i.e., $x_{i,j}, \forall (i,j) \in \mathcal{A} $, to minimize the total flow cost while satisfying the node flow balance and the arc capacity constraint. Given the amount of flow on each arc, the maximal flow of a node $i$ is defined as:
\begin{equation}
	x^f_i = \sum_{j \in \mathcal{O}(i)} x_{i j}+\sum_{j \in \mathcal{I}(i)} x_{j i}, \forall i \in \mathcal{N}.
\end{equation}
Given an arc $(i, j)$, the flow ratio of the arc is defined as:
\begin{equation}
f_{ij} = \max\{\frac{x_{ij}}{x^f_i},\frac{x_{ij}}{x^f_j}\}, \forall(i, j) \in \mathcal{A}.
\end{equation}

In~\cite{ge2021interior}, the authors provide an important property of the MCF problem, which can help to select the basis based on an approximate solution  of~\eqref{MCF_problem}. The property claims that the arc with a larger value of flow ratio is more likely to be included in the basis. Therefore, given an approximate solution, we can sort all arcs according to their flow ratio values and select the first $|\mathcal{N}|$ arcs to create a set as $\{s_1, s_2, \dots, s_{|\mathcal{N}|}\}$, where $|\mathcal{N}|$ is the number of elements in the set $\mathcal{N}$. Intuitively, this set gives the potential arcs which should be included into the basis based on the current approximate solution. 

With this potential set, a column generation basis identification process is used to obtain the feasible basis of the original problem. First of all, the original problem should be converted into an equivalent problem with the standard form. Then, the artificial variable is added based on the big-$M$ method mentioned before. After that, at each iteration, an index set $D_k$ is constructed as: 
\begin{equation}
{D}_{k}  \leftarrow{D}_{k-1} \cup D_{\{s_{1}, s_{2}, \ldots, s_{n_{k}}\}},
%\\\cap\left\{n<i \leq m+n: x_{i}^{k} \text { is nonbasic }\right\}^{c} .
\end{equation}
where $n_k$ is a monotonously increasing sequence of integers. The set $D_{\{s_{1}, s_{2}, \ldots, s_{n_{k}}\}}$ includes the indices of the first $n_k$ arcs in the potential basis set. With the column generation method, the problem at each iteration becomes the following form:
\begin{equation}
\label{column generation}
\begin{aligned}
\min_{x} \quad 	&\sum_{i \in {D_k}} c_{j} x_{j} \\
{\rm s.t.\quad} 	&\sum_{i \in {D_k}} A_{\bullet j} x_{j}=b\\
&x_{j\notin D_k} = 0,\\
&x \geq 0.
\end{aligned}
\end{equation}
%where $A_{\bullet j}$ is the $j$-th column of the $A$.
Once there is no artificial variable in the basis at a particular iteration, the feasible basis for the original problem near the approximate solution is obtained.

This method can be further implemented to obtain an $\epsilon$-optimal feasible solution by changing the update rule of the index set $D_k$ as follows:
\begin{equation}
\label{near optimal update rule}
{D}_{k}  \leftarrow{D
}_{k-1} \cup D_{\{s_{1}, s_{2}, \ldots, s_{n_{k}}\}} \cup \{j: \overline{c}_j < -\epsilon\},
\end{equation}
where $\overline{c}_i$ is the reduced cost.

\emph{Spanning Tree Method}:
Given an MCF problem, the basic feasible solution of this problem should be a feasible tree solution of the graph's spanning tree. Therefore, in the spanning tree method, the task of finding the basic feasible solution can be replaced by the task of constructing a spanning tree with the largest sum of flow ratios. However, if the approximate solution is not accurate, the tree solution may be infeasible. Thus, one more step to take under this condition is to push the infeasible tree solution to a feasible one. For the OT problem, this step can be easily implemented.

Considering the OT problem, the nodes in a directed graph can be divided into the supplier set $\mathcal{S}$ and the consumer set $\mathcal{C}$. Each node in $\mathcal{S}$ has arcs pointing to the nodes in $\mathcal{C}$. The target of the OT problem is to design the amount of flow on each arc to minimize the total cost under given constraints. The general form of the OT problem is given by:
\begin{equation}
\label{ot problem}
\begin{aligned}
\min_{x} \quad 	&\sum_{(i,j) \in \mathcal{S} \times \mathcal{C}} c_{ij} x_{ij} \\
{\rm s.t.\quad} 	&\sum_{j \in \mathcal{C}} x_{i,j}=s_i\\
&\sum_{i \in \mathcal{S}} x_{i,j}=d_j\\
&x_{i,j} \geq 0, \forall (i,j) \in \mathcal{S} \times \mathcal{C}
\end{aligned},
\end{equation}
where $c_{i,j}$ is the cost on the arc $(i,j)$. The variables $s_i$ represents the supply value of a given node $i \in \mathcal{S}$, while $d_j$ is the demand value of a node $j \in \mathcal{C}$.

For any OT problem, if the amount of flow on a given arc $(i,j)\in \mathcal{S} \times \mathcal{C}$ is negative or infeasible, there always exists a new flow design that can be feasible. The detailed method for getting the new feasible design is given in~\cite{ge2021interior} and is omitted here. After this step, the spanning tree method can be used to construct the basic feasible solution as mentioned before.

\section{Initialization in Dual Simplex}
\label{chap:dual_initial}
Other than finding a primal feasible basis, we require a dual feasible point to initialize the dual simplex algorithm. The initialization method in dual simplex can be classified into two types. Methods of the first type solve either the original problem or an auxiliary one with the conventional simplex method. The obtained solution is directly a feasible basic solution to the dual problem. On the other hand, methods of the second type use a modified simplex method instead. Extra operations are implemented at each iteration of the conventional simplex method.

\subsection{Generate a Dual Feasible Basis with Simplex}
\subsubsection{The Most-Obtuse-Angle Row Rule.}\label{chap:obtuse}
The most-obtuse-angle rule was first proposed by Pan in \cite{pan1990practical} and was later analyzed in \cite{pan1994ratio,pan1997most} for its application in achieving dual feasibility. It is a dual version of the most-obtuse-angle column rule. It is inspired by the same observation illustrated in Figure \ref{fig:observation}. In each iteration, one dual infeasible variable is moved from the non-basis into basis, while the leaving variable is selected according to the angle with the down-hill edge determined by the entering variable. If the down-hill direction is close to the negative direction of the objective function, the newly constructed basis is more likely to be an optimal basis. The detailed process is shown in the following:
\begin{enumerate}[1)]
\item Select the entering index $q=\arg\min_{j\in N}s_j$. If $s_q\geq0$, the basis is already dual feasible and go to step \ref{step:final2}).\label{step:initial1}
%\item Compute the transformed pivot column $\Delta x_B=A_B^{-1}A_Ne_q$. If $\Delta x_B\leq0$, the algorithm terminates with (primal) unboundedness. Otherwise, select the leaving index $p=\arg\max_{i\in B}\Delta x_{i}$.
\item Compute $\Delta x_B=A_B^{-1}A_Ne_q$. If $\Delta x_B\leq0$, the algorithm terminates with (primal) unboundedness. Otherwise, select the leaving index $p=\arg\max_{i\in B}\Delta x_{i}$.
\item Perform pivoting $B \gets B\cup\{q\}\backslash\{p\}$ and go to step \ref{step:initial1}).
\item Apply dual simplex to compute the optimum. \label{step:final2}
\end{enumerate}

Though the feasibility of other variables cannot be preserved and cycling may arise during iterations, due to the geometrical benefit of this method, the computational efficiency has been confirmed in \cite{pan1997most,koberstein2007progress}.

\subsubsection{Right-Hand Side Modifications.}
This method was proposed by \cite{pan1999new}, and can be regarded as the dual version of the cost modification method, thus having a similar basic idea. By modifying the right-hand side (RHS) coefficient $\bar{b}\coloneqq A_B^{-1}b$, or equivalently $b$, of the initial problem, the obtained one becomes primal feasible. After solving this modified problem and recovering the right-hand side data, a dual feasible basis can be derived.

Recall $I$ defined in \ref{chap:primal infeasible} which denotes the set of primal infeasible variables. The modified RHS takes
\begin{equation}
	\hat{\bar{b}}_i=\left\{
	\begin{aligned}
		&\delta_i,\quad  \text{if }i \in I \\
		&\bar{b}_i, \quad \text{otherwise,}
	\end{aligned}
	\right.
\end{equation}
where $\delta_i$ is a small positive number to avoid degeneracy. Since we have $x_{i\in I}=\delta_i\geq 0$ for the modified problem, the basis now becomes primal feasible. After solving the modified problem, i.e., $s_N\geq 0$, the original problem can be recovered by restoring the value of $b$. The recovered problem is clearly dual feasible.
\subsection{Generate a Dual Feasible Basis with Modified Simplex}
%\subsubsection{Big-M in Dual}
\subsubsection{Dual Infeasibility-Sum Method.}\label{chap:infeasible-sum}
Though the main idea of dual infeasibility-sum method is not new \citep{maros1986general}, its efficiency was re-analyzed in \cite{maros2003piecewise,koberstein2007progress}. Starting with a dual infeasible basis, this method aims to generate a feasible one by minimizing the infeasibility-sum. Different from the most-obtuse-angle row rule, this method can guarantee a monotonous decrease in infeasibility-sum during iterations.  

Recall the dual basic solution in \eqref{eq:dual-basis}. The basis is said to be infeasible if there exists $s_{j\in N}<0$. The infeasibility-sum method intends to maximize the second term in the summation of such variables, i.e.,
\begin{equation}\label{eq:obj}
\sum_{j\in J}s_j = \sum_{j\in J}c_j-(\sum_{j\in J}A_{\bullet j})^Ty,
\end{equation}
where the objective function is called infeasibility-sum and $J$ is the set of dual infeasible variable defined in \ref{chap:cost-mod}. Therefore, we can construct an auxiliary dual problem based on \eqref{eq:dual-basis} as follows:
\begin{equation}\label{eq:min_infeasible}
\begin{aligned}
\max_{y, s_B, s_N} \quad 	&-(\sum_{j\in J}A_{\bullet j})^Ty\\
{\rm s.t.\quad} 	&A_B^Ty+s_B=c_B\\
						&A_N^Ty+s_N=c_N\\
               			&s_B,s_{N/J}\geq0
\end{aligned}.
\end{equation}
Note that \eqref{eq:min_infeasible} is an LP problem as well, and its basis is exactly the same as that of the original dual problem \eqref{eq:dual}. Nevertheless, as the constraint $s_N\geq0$ becomes $s_{N/J}\geq0$, the basis is dual feasible. Therefore, the dual simplex method can be applied to generate an dual feasible basis of the original problem. One thing needs to mention is that as the constraints are slightly different here, the selection of entering index is slightly modified as well.
%\footnote{As the constraints are slightly different here, the selection of entering index is slightly modified.}.
%
\begin{thm}
Let $\{\bar{x}_B,\bar{x}_N\}$ be the primal basic solution of the auxiliary problem, i.e., $\bar{x}_B=-A_B^{-1}\sum_{j\in J}A_{\bullet j}$ and $\bar{x}_N=0$. The original problem is (primal) unbounded if $\bar{x}_B\geq0$.
\end{thm}

With the above theorem, the dual infeasibility-sum method proceeds in the following steps:
\begin{enumerate}[1)]
\item If $s_N\geq0$, the basis is dual feasible and go to step \ref{final step2}). \label{first step2}
\item Form an auxiliary problem with respect to the current basis and compute $\bar{x}_B=-A_B^{-1}\sum_{j\in J}A_{\bullet j}$. If $\bar{x}_B\geq0$, stop with primal unboundedness. Otherwise, apply one iteration of the modified dual simplex method and go to step \ref{first step2}).
\item Apply dual simplex to compute the optimum of the original problem.\label{final step2}
\end{enumerate} 
\subsubsection{Artificial Bounds.}
The artificial bounds method is a dual version of the big-$M$ method. Simply put, this method adds extra upper bounds to the non-basic variables so that the dual problem has a straightforward feasible basis.

Consider the following problem where there exist newly added artificial bounds for non-basic variables compared with \eqref{eq:primal}:
\begin{equation}\label{eq:artificial-bounds}
\begin{aligned}
\min_{x_B, x_N} \quad 	&c_B^Tx_B+c_N^Tx_N\\
{\rm s.t.\quad}	&A_Bx_B+A_Nx_N=b\\
					&x_B,x_N\geq0\\
					&x_N\leq M
\end{aligned}
\end{equation}
where $M\gg0$ is a large number. Since the artificial bounds can be rewritten as a constraint, i.e., $x_N+x_a=M$ with $x_a\geq0$, where $x_a\in\mathbb{R}^{n-m}$ is the introduced artificial variable, it can be easily verified that the associated dual problem of \eqref{eq:artificial-bounds} is:
\begin{equation}\label{eq:dual-big-m}
\begin{aligned}
\max_{y, y_a, s_B, s_N} \quad		&b^Ty+M{\bf 1}^T_{n-m}y_a\\
{\rm s.t.\quad} 	&A_N^Ty+y_a+s_N=c_N\\
					&A_B^Ty+s_B=c_B\\
					&s_B,s_N\geq0\\
					&y_a\leq0
\end{aligned},
\end{equation}
where $y_a\in\mathbb{R}^{n-m}$ is an artificial variable, serving as the counterpart of $x_a$. Recall the basic solution of the original dual problem where $s_N=c_N-A_N^Ty$. We note that such solution is infeasible if there exists $c_j-A_{\bullet j}^Ty<0$ for any $j\in N$. Replacing such basic variable with $y_a$ in \eqref{eq:dual-big-m}, we can directly obtain a feasible basis of \eqref{eq:dual}.

\section{Conclusion and Future Work}
\label{chap:future}
In this survey, we investigate and summarize existing works related to the simplex initialization from the aspects of the primal simplex and the dual simplex, respectively. A comparison about the discussed methods for generating the initial or staring point is summarized in Table 1. We also discussed some methods to accelerate the initialization in some specific conditions.

\begin{table}[h]
\tiny
\centering
\caption{Comparison of initialization methods in simplex.}
\begin{tabular}{|c|c|c|c|c|c|c|c|}
\hline
\multicolumn{2}{|c|}{Methods}                                      & Sparsity & Triangularity & Creation Time & Feasibility & Optimality    \\ \hline
\multicolumn{2}{|c|}{Two Phase}                                    & $\times$ & $\times$      & long          & $\surd$     & $\times$      \\ \hline
\multicolumn{2}{|c|}{Big-$M$}                                      & $\times$ & $\times$      & long          & $\surd$     & $\times$      \\ \hline
\multicolumn{2}{|c|}{Nonfeasible Basis Method\citep{nabli2009overview}}& $\times$ & $\times$  &       long        & $\surd$     &        $\times$        \\ \hline
\multicolumn{2}{|c|}{Algebraic Initialization \citep{nabli2015algebraic}}& $\times$ &$\surd$  &        long       & $\times$    & $\times$      \\ \hline
\multirow{2}{*}{Modifications}     & Cost \citep{wunderling1996paralleler,pan2000primal}&\multirow{2}{*}{$\times$}&\multirow{2}{*}{$\times$}&\multirow{2}{*}{long}&\multirow{2}{*}{$\surd$}&\multirow{2}{*}{$\times$}\\ \cline{2-2} 
                                   & Right-Hand Side \citep{pan1999new}&      &               &               &             &               \\ \hline
\multirow{2}{*}{Most-Obtuse-Angle} & Column Rule \citep{pan1994variant}&\multirow{2}{*}{$\times$}&\multirow{2}{*}{$\times$}&\multirow{2}{*}{long}&\multirow{2}{*}{$\surd$}&\multirow{2}{*}{$\surd$}\\ \cline{2-2} 
                                   & Row Rule \citep{pan1990practical,pan1994ratio,pan1997most}                     &          &               &               &             &               \\ \hline
\multirow{2}{*}{Infeasibility-Sum} & Primal \citep{ping2014linear}&\multirow{2}{*}{$\times$}&\multirow{2}{*}{$\times$}&\multirow{2}{*}{long}&\multirow{2}{*}{$\surd$}&\multirow{2}{*}{$\times$}\\ \cline{2-2} 
                                   & Dual \citep{maros2003piecewise,koberstein2007progress}                          &          &               &               &             &               \\ \hline
\multicolumn{2}{|c|}{Logical Basis \citep{chvatal1983linear}}      & $\surd$  & $\surd$       & short         & $\times$    & $\times$      \\ \hline
\multicolumn{2}{|c|}{Crash Basis \citep{maros1998strategies}}      & $\times$ & $\surd$       &        short       & LTSF        & $\times$      \\ \hline
\multicolumn{2}{|c|}{CPLEX Basis \citep{bixby1992implementing}}    & $\surd$  & $\times$      &       short        & $\surd$     &       $\times$         \\ \hline
\multicolumn{2}{|c|}{Tearing Basis \citep{gould1989new}}           & $\times$ & $\times$      &        long       & $\surd$     & $\times$      \\ \hline
\multicolumn{2}{|c|}{Cosine Criterion \citep{junior2005improved,hu2007note}}&$\times$&$\times$& long          & $\times$    & $\surd$       \\ \hline
\multicolumn{2}{|c|}{Triangular and Fill-Reducing Basis \citep{ploskas2020triangulation}}&$\surd$&$\surd$&short&$\times$    & $\times$      \\ \hline
\multicolumn{2}{|c|}{  Idiot Crash Algorithm \citep{galabova2020idiot}}&$\times$&$\times$&short&$\times$    & $\surd$      \\ \hline
\multicolumn{2}{|c|}{ $\epsilon$-Optimality Search Direction algorithm \citep{luh2002efficient}}&$\times$&$\times$&long &$\surd$    & $\surd$      \\ \hline
\multicolumn{2}{|c|}{ Hybrid-LP \citep{al2011hybrid}}&$\times$&$\times$&long&$\surd$& $\surd$      \\ \hline
\end{tabular}
\end{table}

For the future work, there are mainly two directions. One is based on the conventional initialization methods introduced in this survey, which include: investigate the efficiency of the existing algorithms for dealing with large scale LPs; provide some generalized algorithms; discover implementation techniques to further improve the efficiency of existing algorithms; formulate the dual version of algorithms provided in the primal part, and investigate their performance. The other one is based on advanced learning technologies. As machine learning techniques are gaining more and more research and development, many fields are exploring how to combine existing methods with advanced machine learning techniques to get improved methods. However, there is almost no research in the field of simplex initialization that utilizes learning-based methods to improve the efficiency of solving LPs. Therefore, it is important to investigate and fill the gap in this field. In the following part, we propose several potential directions for utilizing the learning-based methods to improve the simplex initialization.

\subsection{Learning-Based Initial Basis Construction}
In the last two sections, we have discussed many approaches to construct the initial basis or point from different perspectives. For example, the logical basis is one of the easiest ways to obtain a feasible initial basis, since it directly uses all logical variables to construct a basis. Different from the logical basis, the cosine criterion method obtains the initial basis based on the constraints with the minimum angle 
to the objective function. 
%, without considering the sparsity or triangularity of the basis. 
Other methods such as triangular and fill-reducing basis focus on constructing a triangular and sparse basis. Choosing the appropriate initialization methods for each LP based on its features is one of the possible ways to improve the overall efficiency. With advanced machine learning technologies, we can design a classifier for automatically selecting the appropriate simplex initialization methods to construct the initial basis. In addition, we can design a classifier for determining whether a variable should be included into the basis directly.

% In this subsection, we will present two learning-based methods. One is for determining whether a variable should be included into the basis. The other is for automatically selecting the appropriate simplex initialization methods to construct the initial basis. These methods are based on some machine learning and deep learning classification methods. 

\subsubsection{Feature Design.}\label{sec:FE-GE}
For different LPs, we first need to design some features to distinguish them. There are two types of features, namely, self-designed features and graph embedding features.

\emph{Self-Designed Features}: Consider the LP standard form given in~\eqref{eq:primal}, the self-designed features can be further divided into the problem-dependent and the problem-independent features. The dimension of the problem-dependent feature changes when given different LPs, while the dimension of the problem-independent feature does not change. The details of the features can be designed as follows.

\begin{enumerate}[1)]
\item Problem-dependent features: the matrix $A$, the vectors $c, b$, the dimensions $m, n$.
\item Problem-independent features: the sparsity of the matrix $A$, which is quantified by the percentage of zeros in $A$; the sparsity of the vector $b$, which is defined by the percentage of zeros in $b$; the triangularity of the matrix $A$, which is designed as:
\begin{equation}
\label{triangularity}
\frac{|z_U - z_L|}{\max(z_U, z_L)},
\end{equation}
where $z_U$ and $z_L$ are the percentages of non-zero elements in the upper part and the lower part of the matrix $A$, respectively.
\end{enumerate}

\emph{Graph Embedding Features:}
%\subsection{Feature extraction based on graph embedding}
Recent advancements in graph embedding have drawn the attention of researchers in various fields. Graph embedding provides a new approach to obtain low-dimensional features of nodes in a graph, which is significant for reducing computational complexity. Intuitively, the main idea of graph embedding is to find a mapping from a node in a graph to a low-dimensional feature representation. This method has been applied in diverse areas such as social sciences, linguistics, biology, etc. Recently, graph embedding has been used to improve the performance of the branch-and-bound algorithm for solving mixed integer programming (MIP)~\citep{nair2020solving}.

\begin{figure}
	[htbp]
	\centering
	\includegraphics[width=5in]{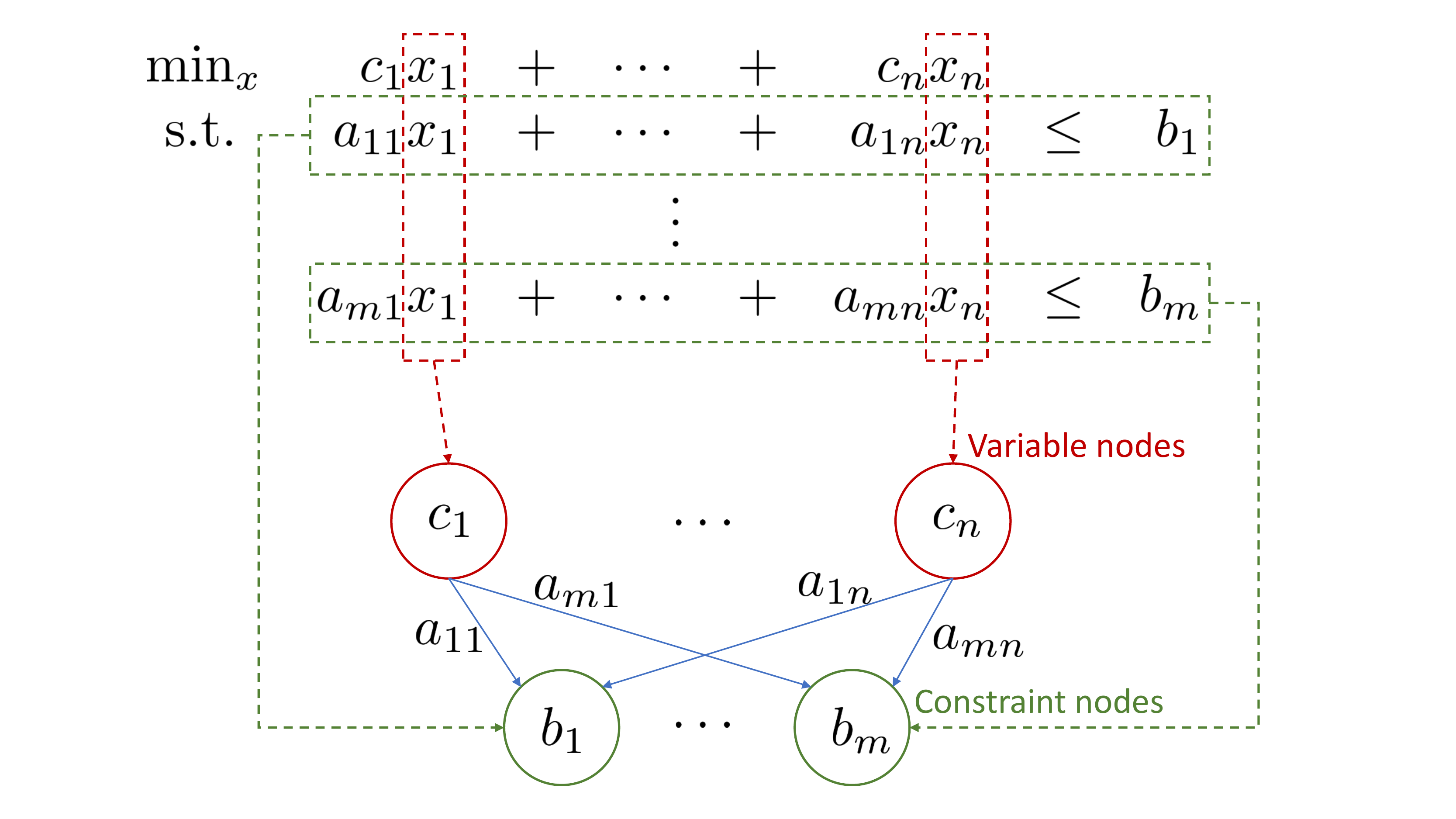}
	\caption{Bipartite graph representation of LP}
	\label{fig:bipartite}
\end{figure}

Motivated by~\cite{selsam2018learning,gasse2019exact,ding2020accelerating} and~\cite{nair2020solving}, we can also represent an LP as a bipartite graph. Considering the standard form \eqref{eq:primal} of LP, the following bipartite graph (Figure \ref{fig:bipartite}) can be constructed.
In the graph, one partition has $n$ (variable) nodes, which represent the $n$ variables to be optimized, and the other has $m$ (constraint) nodes, which represent the $m$ constraints in the standard form of LP. If a variable appears in a constraint, there will exist an edge between the corresponding variable node and constraint node, and the edge is weighted by the corresponding entries of the matrix $A$. The objective coefficients $\{c_1, \ldots, c_n\}$, the right-hand side of the constraints $\{b_1, \ldots, b_m\}$, and the non-zero entries of the matrix $A$ can be utilized as scalar ``features'' of the variable nodes, the constraint nodes, and the edges, respectively. 

The bipartite graph representation of LP enables the feature extraction through graph embedding. Generally, graph embedding methods can be devided into two types: homogeneous graph embedding (e.g., DeepWalk~\citep{perozzi2014deepwalk}, LINE~\citep{tang2015line}, Node2vec~\citep{grover2016node2vec}, and VGAE~\citep{kipf2016variational}) and heterogeneous graph embedding (e.g., Metapath2vec~\citep{dong2017metapath2vec} and DMGI~\citep{park2020unsupervised}). 
Since the nodes in the constructed graph have two different types (variable nodes and constraint nodes), the graph should be considered as a heterogeneous graph. Therefore, all the existing heterogeneous graph embedding methods are possible candidates for learning the low-dimensional features of the nodes. 
However, most methods are not specific to bipartite graphs, and thus may not capture their unique structural characteristics. Some methods (e.g., IGE~\citep{zhang2017learning}, PinSage~\citep{ying2018graph}, BiNE~\citep{gao2018bine}, FOBE~\citep{sybrandt2019fobe}, BiGI\citep{cao2021bipartite}) are specially designed for bipartite graphs, and they may be better choices for our problem.

%\subsubsection{Deep Learning based Graph Embedding}
One subfield of graph embedding is deep learning (DL) based graph embedding~\citep{cai2018comprehensive, goyal2018graph}. 
Based on whether random walk is adopted to sample paths from a graph, DL based methods can be divided into two categories. One is DL based graph embedding with random walk, including DeepWalk~\citep{perozzi2014deepwalk}, Node2vec~\citep{grover2016node2vec}, etc. The other is DL based graph embedding without random walk, including SDNE~\citep{wang2016structural}, DNGR~\citep{cao2016deep}, graph convolutional networks (GCN)~\citep{kipf2016semi}, VGAE~\citep{kipf2016variational}, etc. Compared with factorization-based embedding (e.g., graph factorization (GF)~\citep{ahmed2013distributed}, GraRep~\citep{cao2015grarep}, and HOPE~\citep{ou2016asymmetric}), DL based graph embedding has wider application due to its robustness and effectiveness. Since the corresponding nodes of variables appearing in the same constraint are two-hop neighbors to each other in the constructed bipartite graph, the second (or higher) order proximity is expected to be preserved. Therefore, DL based methods will be preferred for our problem due to their good performance in preserving higher-order properties of graphs.

\subsubsection{Initial Basis Construction with Deep Learning Based Classification.}
\begin{figure}
	[htbp]
	\centering
	\includegraphics[width=5in]{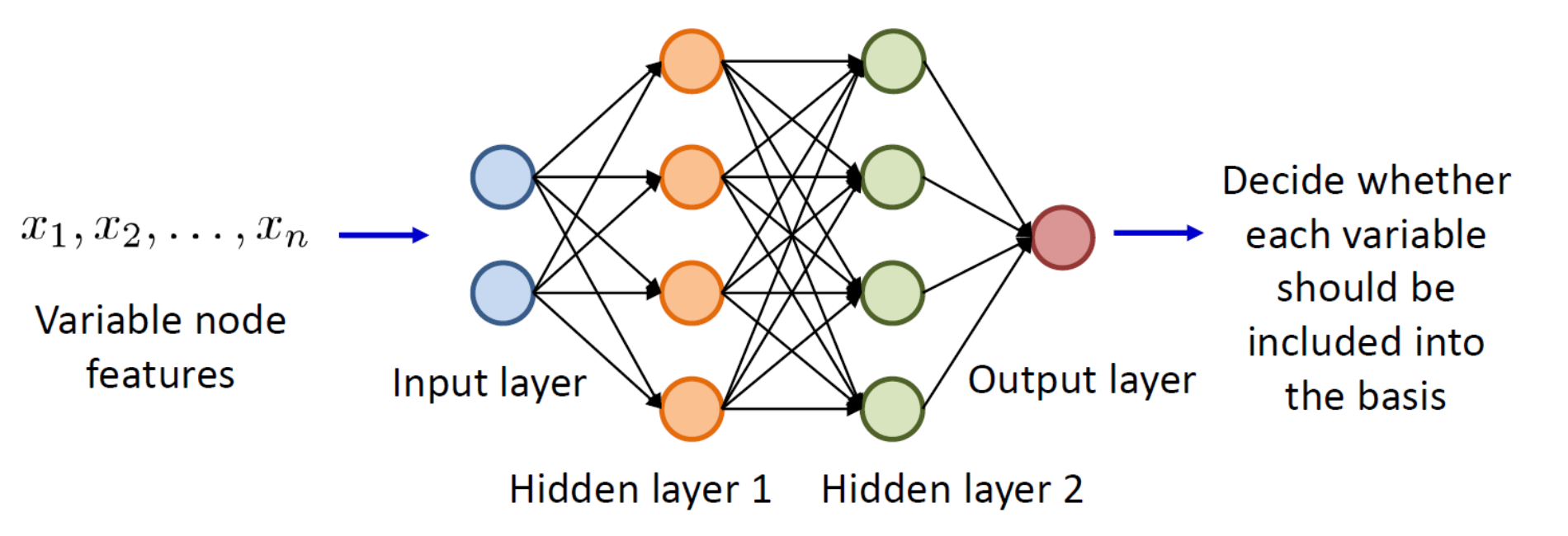}
	\caption{Classification of initial basic variable selection}
	\label{fig:class_initial_variable_select}
\end{figure}

The main purpose of this subsection is to provide a classification mechanism based on a deep neural network, which can divide variables into basic variables and non-basic variables. The input of the neural network is the feature of each variable node, 
%, i.e., a low-dimensional vector $\boldsymbol{v}_i\in\mathbb{R}^d, i\in\{1,\ldots,n\}$ 
which can be obtained through graph embedding. The output is the probability that the corresponding variable should be selected as a basic variable. The architecture is given in Figure~\ref{fig:class_initial_variable_select}.

To train such a neural network, enough training data pairs are required. Each training data pair consists of the feature of a variable node and the label indicating whether its corresponding variable is a basic variable, i.e., the label is ``1'' if it is a basic variable, otherwise ``0''. As we have mentioned before, the features can be obtained by graph embedding. However, how to obtain the labels of variable nodes can be a tricky problem. One approach to obtain the labels is to exactly solve the LP problem, then the type (basic/non-basic) of the variable when reaching the optimality can serve as its label. Obviously, this approach is very costly. Another approach is to obtain the label by the initialization methods introduced before, but the obtained label can be inappropriate since some of these methods cannot guarantee the feasibility or the nonsingularity of the derived basis. If enough training data have been obtained, the feature and the label can serve as the input and the output of the deep neural network, respectively. The architecture of the deep neural network, including the number of layers, the activation function, the loss function, etc., should be further designed. 

The trained neural network can be used to select basic variables for LPs. The main steps are as follows: first, construct the graph representation for the LP to be solved; second, learn the low-dimensional features of the variable nodes via graph embedding; third, input the derived features into the neural network to obtain probability outputs, then sort all variables in descending order of their corresponding outputs and select the first $m$ variables as basic variables.

\subsubsection{Initialization Method Selection with Deep Learning Based Classification.}
\label{sub_classification_selectmethod}
\begin{figure}
	[htbp]
	\centering
	\includegraphics[width=5in]{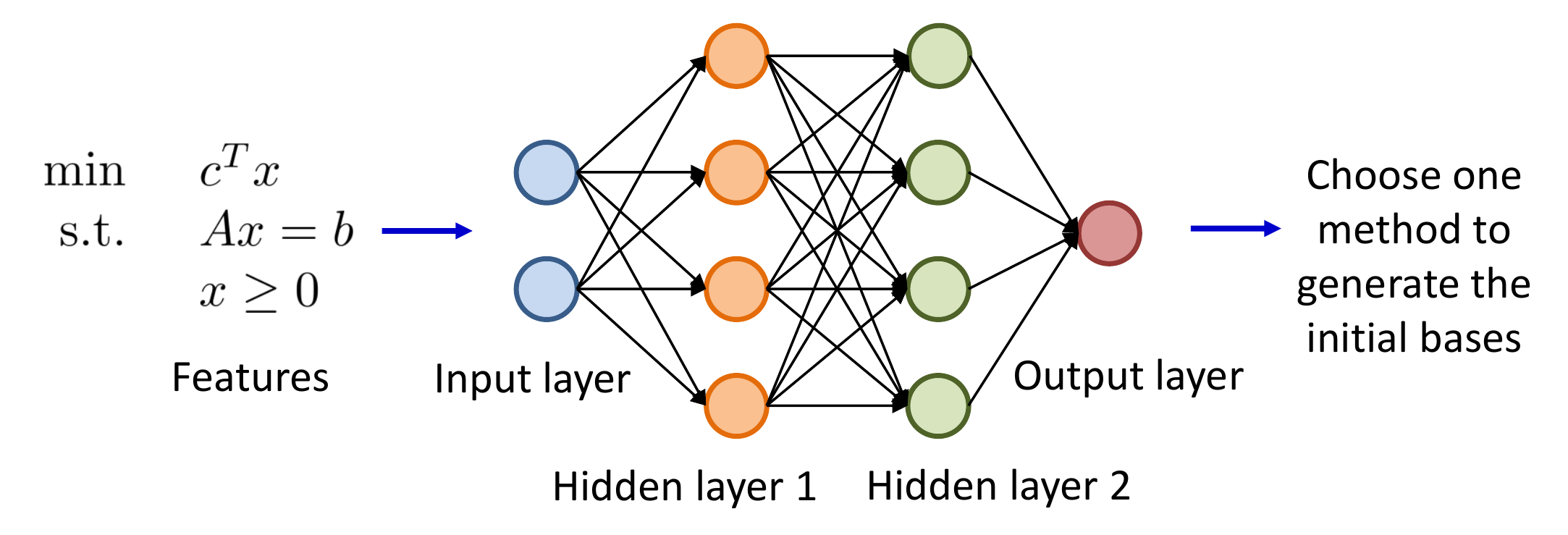}
	\caption{Classification of initialization method selection}
	\label{fig:class_initial_method_select}
\end{figure}

The main idea of the learning-based classification for the initialization method selection is illustrated in Figure~\ref{fig:class_initial_method_select}.

The input of the classifier is the designed features of different LPs. These features can be the self-designed features, the graph embedding features, or a combination of different features, depending on the classification performance. The output of the classifier is a probabilistic distribution of choosing different simplex initialization methods to construct the initial basis. Given an LP, the method with the highest probability will be selected as the optimal initialization method. 

To train the classifier, we first collect the total computation time with different candidate methods and label the method with the shortest time as ``1''. The others can be labelled as ``0''. Then these labels are regarded as the training dataset to train the classifier.

% We propose to consider two groups of candidate initialization methods as follows:
% \begin{enumerate}[1)]
% \item Logical basis, crash basis, CPLEX basis, and the tearing algorithm: These basis construction methods are all applicable to LPs with the same form. However, few works have been done to investigate the optimal selection among these initialization methods. With the proposed classifier, we can automatically select the optimal method to do the initialization. 

% \item Logical basis and the triangular and fill-reducing basis: The logical basis is one of the simplest methods to construct the basis. It construct the initial basis with the artificial variables directly. The triangular and fill-reducing basis, on the other hand, focuses on constructing a sparse and triangular basis. For different LPs with different features or structures, the learning-based classifier can help to choose the appropriate basis.

% \end{enumerate}

\subsection{Learning-Based Starting Point Construction}
In the Simplex Initialization section, we have covered some methods for finding an improved starting point for the simplex method. These methods design different interior directions to do the iteration and obtain the improved point. However, the trade-off between the iteration steps and the total computation time has not been well studied. In this subsection, we will propose the reinforcement learning (RL)-based method for investigating the trade-off between the iteration steps and the total computation time in methods for finding the improved starting point. 

\subsubsection{State Design.}
In the RL, an action is selected based on the current state. After the action selection, a reward and the next state will be returned. This process is then repeated in the following steps. In our problem, the action selection is equivalent to designing whether the iteration will continue or stop. The state includes two parts. One is the feature of the given LP. The feature can be the self-designed one or the graph embedding one. The other part is the state related to the current improved point. For example, in the $\epsilon$-optimality search direction method, the iteration direction is designed based on the active constraint of the current improved point. Therefore, the active constraint can be included in the state. The search direction can also be considered as part of the state. Based on the designed state, the action selection is learned by RL algorithms.

\subsubsection{Starting Point Construction with RL.}
\begin{figure}
	[htbp]
	\centering
	\includegraphics[width=5in]{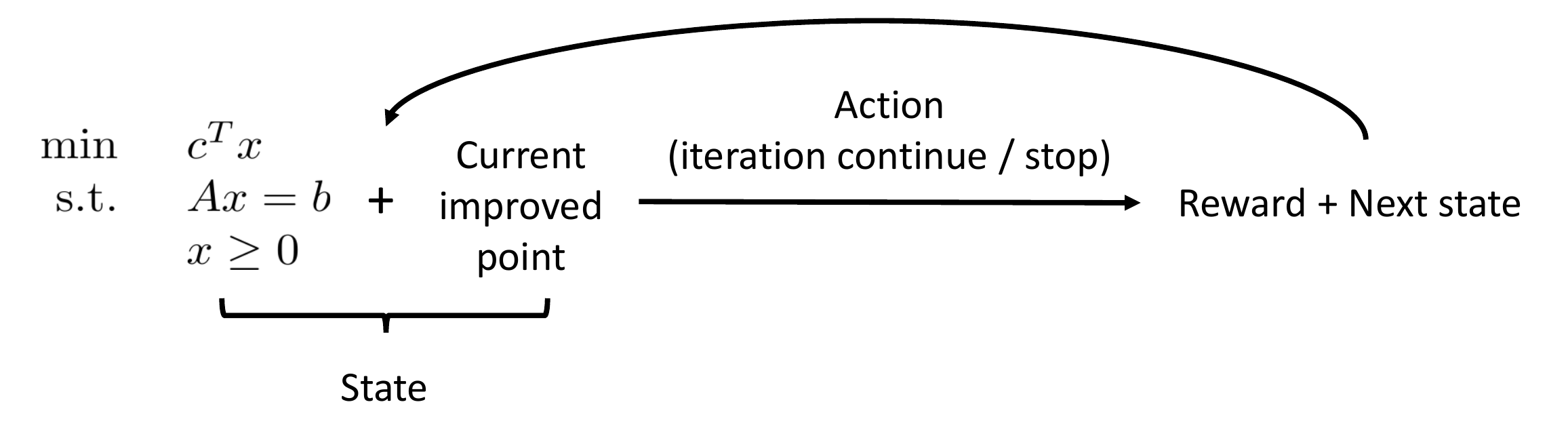}
	\caption{Learning-Based Starting Point Construction}
	\label{fig:rl_future_work}
\end{figure}

The general process of the learning-based starting point construction is given in Figure~\ref{fig:rl_future_work}. At each step, based on the state and the selected action, a reward is returned. The main goal is to maximize the cumulative reward. The reward can be designed based on the computation time of each iteration.

According to the dimensionality of the designed state, different RL methods can be chosen.  When the state dimension is small, we can choose a basic learning-based method such as Q-learning algorithm. When the state space is large or continuous, we can choose some advanced algorithms like the deep Q-network algorithm or double deep Q-network algorithm, etc.

\subsection{Learning Heuristics in Simplex Initialization}

Most of methods we mentioned before have included some heuristic rules or heuristic parameters. In these methods, the selection of these rules and parameters are almost done manually. To improve the efficiency of different simplex initialization methods, we can try to improve the heuristic designs from two perspectives as follows.
\begin{enumerate}[1)]
\item Select the optimal heuristic rules among candidate rules: When there exist several different heuristic rules, we can construct a classification method to choose the better one. Specifically, we can label these candidate rules based on their computation time.

\item Generate the heuristic parameters automatically: When there are some LPs with appropriate heuristic parameters, we can formulate a regression problem to obtain a mapping from the problem features to the heuristic parameters. Then when a new LP is given, we can generate the suitable heuristic parameter automatically based on the regression model.
\end{enumerate}

%%%%%%%%%%%%%%%%%%%%%%

\bibliographystyle{unsrtnat}
\bibliography{references}  %%% Uncomment this line and comment out the ``thebibliography'' section below to use the external .bib file (using bibtex) .





% \end{thebibliography}

\end{document}